\documentclass[review,hidelinks,onefignum,onetabnum]{siamart250106}

\nolinenumbers

\usepackage{amsmath,amsfonts,upgreek}

\newsiamremark{ass}{Assumption}
\newsiamremark{ex}{Example}
\newsiamremark{rem}{Remark}

\newcommand{\be}{\begin{equation*}}
\newcommand{\ben}{\begin{equation}}
\newcommand{\ee}{\end{equation*}}
\newcommand{\een}{\end{equation}}

\newcommand{\sgn}{\operatorname{sgn}}

\newcommand{\bbr}{\mathbb{R}}

\newcommand{\EXP}{\operatorname{\textsf{\upshape E}}}
\newcommand{\PR}{\operatorname{\textsf{\upshape P}}}

\newcommand{\acal}{\mathcal{A}}
\newcommand{\fcal}{\mathcal{F}}

\newcommand{\JE}{J^{\mathrm e}}
\newcommand{\JP}{J^{\mathrm p}}

\newcommand{\cc}{{\mathrm c}}

\newcommand{\di}{\mathrm{d}}

\newcommand{\black}{\color{black}}
\newcommand{\blue}{\color{black}}

\begin{document}

\title{Singular stochastic control problems motivated
by the optimal sustainable exploitation of an
ecosystem\thanks{Research supported by the LMS
Grant Scheme 4:\ Research in Pairs (Ref:\ 41920).}}

\author{Gechun Liang\thanks{Department of Statistics,
University of Warwick, Coventry, CV4 7AL, UK,
\email{g.liang@warwick.ac.uk}.}
\and
Zhesheng Liu\thanks{Department of Applied Mathematics,
The Hong Kong Polytechnic University,  Hung Hom,
Kowloon, Hong Kong, China, \email{liuzhesheng@outlook.com}.}
\and
Mihail Zervos\thanks{Department of Mathematics,
London School of Economics, Houghton Street, London
WC2A 2AE, UK, \email{mihalis.zervos@gmail.com}.}
}

\maketitle

\begin{abstract}
We derive the explicit solutions to singular stochastic control
problems of the monotone follower type with
(a) an expected discounted criterion,
(b) an expected ergodic criterion and
(c) a pathwise ergodic criterion.
These problems have been motivated by the optimal sustainable
exploitation of an ecosystem, such as a natural fishery.
Under general assumptions on the diffusion coefficients,
the discounting rate function, the running payoff function
and the marginal profit of control action, we show that the
optimal strategies are of a threshold type.
We solve the three problems by first constructing suitable
solutions to their associated HJB equations,
which take the form of quasi-variational inequalities with
gradient constraints.
In the cases of the ergodic control problems,
we also use a suitable new variational argument.
Furthermore, we establish the convergence
of the solution of the discounted control problem to the
one of the ergodic control problems as the discounting
rate function tends to 0 in an Abelian sense.
\end{abstract}



\headers{Singular stochastic control problems}
{G. Liang, Z. Liu, and M. Zervos}

\section{Introduction}

We consider a stochastic dynamical system with a positive state
process that satisfies the stochastic differential equation
\ben
\di X_t^\zeta = b(X_t^\zeta) \, \di t - \di \zeta_t  + \sigma
(X_t^\zeta) \, \di W_t , \quad X_{0-}^\zeta = x  > 0 , \label{SDE}
\een
where $W$ is a standard one-dimensional Brownian motion
and $\zeta$ is a controlled c\`{a}dl\`{a}g increasing process.
With each controlled process $\zeta$, we associate the
expected discounted performance index
\ben
I_x (\zeta) = \EXP_x \Biggl[ \int _0^\infty e^{-\Lambda_t^\zeta}
h (X_t^\zeta) \, \di t + \int _0^\infty e^{-\Lambda_t^\zeta}
k(X_t^\zeta) \circ \di \zeta_t \biggr] ,  \label{JEd}
\een
the expected long-term average performance index
\ben
\JE_x (\zeta) = \limsup _{T \uparrow \infty} \frac{1}{T}
\EXP \biggl[ \int _0^T h (X_t^\zeta) \, \di t + \int _0^T
k(X_t^\zeta) \circ \di \zeta_t \biggr] ,
\label{JE}
\een
as well as the pathwise long-term average performance
criterion
\ben
\JP_x (\zeta) = \limsup _{T \uparrow \infty} \frac{1}{T}
\biggl( \int _0^T h (X_t^\zeta) \, \di t + \int _0^T k(X_t^\zeta)
\circ \di \zeta_t \biggr) , \label{JP}
\een
where
\begin{gather}
\Lambda _t^\zeta = \int _0^t r(X^\zeta_u) \, \di u
\label{Lam} \\
\text{and} \quad
\int _0^T k(X_t^\zeta) \circ \di \zeta_t = \int _0^T k(X_t^\zeta)
\, \di \zeta_t^\cc + \sum _{0 \leq t \leq T} \int _0^{\Delta \zeta_t}
k(X_{t-}^\zeta - u) \, \di u . \label{kdzeta-int}
\end{gather}
In the last of these definitions, $\zeta^\cc$ is the continuous
part of the c\`{a}dl\`{a}g increasing process $\zeta$ and
$\Delta \zeta_t = \zeta_t - \zeta_{t-}$, with the convention
that $\zeta _{0-} = 0$.
The objective of the resulting singular stochastic control
problems is to maximise each of the objective criteria
(\ref{JEd}), (\ref{JE}) and (\ref{JP}) over all admissible
controlled processes $\zeta$.

The control problems defined by (\ref{SDE})--(\ref{JP}) have
been motivated by the sustainable exploitation of an ecosystem,
such as a forest or a natural fishery.
In such a context, $X$ models the population level process
of a harvested species, while $\zeta_t$ is the total amount
of the species that has been harvested by time $t$.
The function $k > 0$ in (\ref{JEd})--(\ref{JP}) models the profit
made from each unit of the harvested species.
On the other hand, the function $h$ models the utility arising
from having a population level $X_t$ of the harvested species
at time $t$, which could reflect the role that the species plays
in the stability of the overall ecosystem.
Alternatively, the function $h$ can be used to model running
costs.

Motivated by applications to the optimal harvesting of
stochastically fluctuating populations, similar singular
stochastic control problems with $h=0$, constant $k$
and with a discounted performance criterion with constant
$r$ have been studied by
Alvarez~\cite{ALV00, ALV01},
Alvarez and Shepp~\cite{AS}, and
Lungu and {\O}ksendal~\cite{LO97}.
Extensions of these earlier works have been studied by
Framstad~\cite{FRA}, who considers a state process $X$
with jumps,
Song, Stockbridge and Zhu~\cite{SSZ}, who consider a state
process $X$ with regime switching,
Morimoto~\cite{M13}, who considers the finite time horizon
case,
Alvarez, Lungu and {\O}ksendal~\cite{ALO} and
Lungu and {\O}ksendal~\cite{LO01},
who consider multidimensional state processes $X$,
Hening, Tran, Phan and Yin~\cite{HTPY}, who consider
multidimensional state processes $X$ as well as allow for
the modelling of both seeding and harvesting, and
Ga\"{i}gi, Ly\,Vath and Scotti~\cite{GLVS22},
who consider constraints of no-take areas.
On the other hand, control problems
with an expected ergodic performance criterion,
similar to the one that we study here with $h = 0$
and constant $k$, have been solved by
Hening, Nguyen, Ungureanu and Wong~\cite{HNUW},
Alvarez and Hening~\cite{AH22}, as well as
Cohen, Hening and Sun~\cite{CHS22}, who consider
a performance criterion with model ambiguity.
Several other closely related contributions can be
found in the literature of all these papers.

We solve the control problems that we consider by
deriving explicit solutions to their corresponding
HJB equations.
In generalising the special cases arising when $h=0$
and $k$ is constant, our main contributions include
(a) the determination of sufficiently general assumptions
on the functions $h$ and $k$  that give rise to
threshold optimal strategies without making extra
assumptions on the data $b$ and $\sigma$ of the
underlying diffusion,
and (b) the derivation of explicit solutions to the problems'
HJB equations that are way more complicated than the ones
associated with the special case arising when $h=0$
and $k$ is constant
(e.g., we are faced with integral equations, such as
(\ref{Xi}), instead of algebraic equations, such as
the one in Remark~\ref{rem:AH19} from 
Alvarez and Hening~\cite{AH22}).

We derive the solution to the discounted singular stochastic
control problem in Section~\ref{sec:disc-sol} using results
from Liu and Zervos~\cite{LZ24}, who solve the corresponding
discounted impulse control problem.
On the other hand, we solve the ergodic singular stochastic
control problems in Sections~\ref{sec:HJB-sol}
and~\ref{sec:erg-sol}.
The analysis of these problems, which are in the so-called
monotone follower singular stochastic control setting, has
been influenced by
Karatzas~\cite{K83},
Menaldi, Robin and Taksar~\cite{MRT92},
Weerasinghe~\cite{W02}, and
Jack and Zervos~\cite{JZ},
who consider different formulations.
In the solution to the ergodic control problems that we
solve, a notable difficulty arises from the fact that the
solution $(w, \lambda^\star)$ to their corresponding
HJB equation may involve functions
$w$ that are unbounded from below (see
Remark~\ref{rem:why(III)}), which makes the
establishment of a suitable verification theorem
intractable.
We overcome this complication by means of
a variational argument involving suitable
pairs $(w_\lambda, \lambda)$ with bounded
from below functions $w_\lambda$ that converge
to $(w, \lambda^\star)$ as
$\lambda \downarrow \lambda^\star$.
The introduction of this technique is a further
contribution of this paper.

In Section~\ref{sec:Abelian}, we establish the convergence
of the solution to the discounted control problem to the
one of the ergodic control problems as the discounting
rate function $r$ tends to 0 in an Abelian sense.
In particular, we prove that, if $r$ depends on a
parameter $\iota > 0$ and tends to zero as $\iota
\downarrow 0$ in the sense of Assumption~\ref{ass:Ab},
then
\ben
\lim _{\iota \downarrow 0} \beta^\star (\iota)
= \beta^\star , \quad \lim _{\iota \downarrow 0}
r(y; \iota) w(x; \iota) = \lambda^\star
\quad \text{and} \quad
\lim _{\iota \downarrow 0} w'(x; \iota)
= w'(x) \label{Ab-lims0}
\een
for all  $x,y > 0$, where $\beta^\star (\iota)$
(resp., $\beta^\star$) is the threshold point
characterising the optimal strategy of the
discounted problem (resp., the ergodic problems)
and $w (\cdot ; \iota)$ (resp., $(w, \lambda^\star)$)
is the solution to the HJB equation of the discounted
problem (resp., ergodic problems).
In the context of singular stochastic control, Abelian
limits, such as the first two ones in (\ref{Ab-lims0}),
have been obtained for constant $r$ by
Karatzas~\cite{K83},
Weerasinghe~\cite{W07},
Hynd~\cite{H12},
Alvarez and Hening~\cite{AH22}, and
Kunwai, Xi, Yin and Zhu~\cite{KXYZ22}, \blue
as well as by Cao, Dianetti and Ferrari~\cite{CDF23}
in a mean-field game setting,
\black
using different techniques.
To the best of our knowledge, no results
(a) with non-constant discounting rate functions
$r(\cdot; \iota)$, or (b) such as the third limit in
(\ref{Ab-lims0}), exist in the singular stochastic
control literature, \blue
with the exception of
Karatzas~\cite[Proposition~6]{K83}, who
establishes a limit such as the third one
in (\ref{Ab-lims0}) for a model with a standard
Brownian motion and constant $r$.
\black

\section{Problem formulation}
\label{sec:pr-form}

Fix a filtered probability space $\bigl( \Omega, \fcal, (\fcal_t),
\PR \bigr)$ satisfying the usual conditions and carrying a
standard one-dimensional $(\fcal_t)$-Brownian motion $W$.
We consider a biological system, the uncontrolled stochastic
dynamics of which are modelled by the SDE
\ben
\di X_t = b(X_t) \, \di t + \sigma (X_t) \, \di W_t , \quad
X_0 = x  > 0 , \label{SDE-0}
\een
for some deterministic functions $b , \sigma : \mbox{}
]0 ,\infty [ \mbox{} \rightarrow \bbr$.

\begin{ass} \label{A1}
The function $b$ is $C^1$ and the limit $b(0) := \lim
_{x \downarrow 0} b(x)$ exists in $\bbr$.
On the other hand, the function $\sigma$ is $C^1$,
the limit $\sigma (0) := \lim _{x \downarrow 0} \sigma (x)$
exists in $\bbr$ and
\ben
0 < \sigma ^2 (x) \leq C_1(1+x^\eta) \quad \text{for all }
x>0 , \label{sigma-bounds}
\een
for some constant $C_1, \eta >0$.
\end{ass}
This assumption implies that the scale function $p$ and the
speed measure $m$ of the diffusion associated with the SDE
(\ref{SDE-0}), which are given by
\begin{align}
p (1) = 0 \quad \text{and} \quad & p' (x) = \exp
\biggl( - 2 \int _1^x \frac{b(s)} {\sigma ^2(s)} \, \di s \biggr)
, \quad \text{for } x > 0, \label{p} \\
\intertext{and} &m (\di x) = \frac{2}{\sigma ^2(x) p' (x)}
\, \di x , \label{m}
\end{align}
are well-defined.
We also make the following assumption, which, together with
Assumption~\ref{A1}, implies that the SDE (\ref{SDE-0})
has a unique non-explosive strong solution (e.g.,
see Karatzas and Shreve~\cite[Proposition~5.5.22]{KSh}).

\begin{ass} \label{A2}
The scale function $p$ and the speed measure
$m$ defined by (\ref{p}) and (\ref{m}) satisfy
\be
\lim _{x \downarrow 0} p (x) = -\infty , \quad
\lim _{x \uparrow \infty} p (x) = \infty
\quad \text{and} \quad m \bigl( ]0, 1[ \bigr) < \infty .
\ee
\end{ass}
For the solution to the ergodic control problems, we need
the following additional assumption, which implies
that the diffusion associated with the SDE (\ref{SDE-0})
is ergodic.

\begin{ass} \label{A12+}
The integrability condition
\be
\int _0^\infty \bigl( s^\eta + 1 \bigr) \, m (\di s) < \infty
\ee
holds true, where $\eta > 0$ is as in (\ref{sigma-bounds}).
\end{ass}

If the system is subject to harvesting, then its state process
$X$ satisfies the controlled one-dimensional SDE
(\ref{SDE}).

\begin{definition}
An {\em admissible harvesting strategy\/} is any
$(\fcal_t)$-adapted process $\zeta$ with c\`{a}dl\`{a}g
increasing sample paths such that $\zeta _{0-} =0$ and
the SDE (\ref{SDE}) has a unique non-explosive strong
solution.
We denote by $\acal$ the family of all admissible strategies.
\end{definition}

With each admissible harvesting strategy $\zeta \in \acal$,
we associate the expected discounted performance
index $I_x (\zeta)$ given by (\ref{JEd}), the expected
ergodic performance index $\JE_x$ given by (\ref{JE})
and the pathwise performance criterion
$\JP_x$ given by (\ref{JP}).
The objective of the control problem that we consider is to
maximise each of $I_x (\zeta)$, $\JE_x (\zeta)$ and
$\JP_x (\zeta)$ over all $\zeta \in \acal$.

\begin{ass} \label{A3all}
(i)
The function $h$ is $C^1$ as well as bounded from
below and the limit $h(0) := \lim _{x \downarrow 0} h(x)$
exists in $\bbr$.
\smallskip

\noindent (ii)
The function $k$ is positive, bounded and $C^2$.
Also, the limits $k(0) := \lim _{x \downarrow 0} k(x)$
and $k'(0) := \lim _{x\downarrow 0} k'(x)$ both exists in
$\bbr_+$.
\end{ass}

For the discounted control problem, we make
the following additional assumption.

\begin{ass} \label{A4d}
(i)
The discounting rate function $r$ is bounded and $C^1$.
Also, there exists $r_0>0$ such that $r(x) \geq r_0$ for all
$x \geq 0$.
\smallskip

\noindent (ii)
The limit $\lim _{x \downarrow 0} h(x) / r(x)$ exists in
$\bbr$ and
\be
\EXP_x \biggl[ \int _0^\infty e^{-\Lambda_t} \bigl| h(X_t)
\bigr| \, \di t \biggr] < \infty .
\ee

\noindent (iii)
The conditions
\be
\EXP_x \biggl[ \int _0^\infty e^{-\Lambda_t} \bigl\vert
\Theta (X_t) \bigr\vert \, \di t \biggr] < \infty
\quad \text{and} \quad
\limsup _{x \uparrow \infty} \frac{1}{\psi (x)} \int _0^x
k(s) \,  \di s \in \bbr_+
\ee
hold true, where $\psi$ is as at the beginning of
Section~\ref{sec:disc-sol} and
\ben
\Theta (x) = \frac{1}{2} \sigma^2 (x) k'(x) + b(x) k(x)
- r(x) \int _0^x k(s) \,  \di s , \quad \text{for } x > 0 .
\label{Theta}
\een

\noindent (iv)
The limit $\lim _{x \downarrow 0} \Theta (x) / r(x)$ exists
in $\bbr$ and there exists a constant $\xi \in \mbox{}
]0 ,\infty [$ such that
\be
\frac{\di}{\di x} \frac{\Theta (x) +h(x)}{r(x)}
\begin{cases} > 0, & \text{for } x \in \mbox{} ]0, \xi [ ,
\\ < 0 , & \text{for } x \in \mbox{} ]\xi ,\infty[ ,
\end{cases}
\quad \text{and} \quad
\lim _{x \uparrow \infty} \frac{\Theta (x) + h(x)}{r(x)}
< \frac{\Theta (0) + h(0)}{r(0)} ,
\ee
where $\Theta$ is defined by (\ref{Theta}).
\end{ass}

For the ergodic control problems, we make
the following additional assumption.

\begin{ass} \label{A4e}
 (i)
The following integrability condition is satisfied:
\be
\int _0^\infty |h(s)| \, m (\di s) < \infty .
\ee

\noindent (ii)
If we define
\ben
K(x) = \frac{1}{2} \sigma^2 (x) k'(x) + b(x) k(x) , \quad
\text{for } x > 0 , \label{K}
\een
then there exists a constant $\xi \in \mbox{}
]0 ,\infty [$ such that
\be
K'(x) +h'(x) \begin{cases} > 0, & \text{for } x \in \mbox{}
]0, \xi [ , \\ < 0 , & \text{for } x \in \mbox{} ]\xi ,\infty[ ,
\end{cases}
\quad \text{and} \quad
\lim _{x \uparrow \infty} \bigl( K(x) + h(x) \bigr) <
K(0) + h(0) .
\ee
\end{ass}

\begin{rem} \label{rem:Kk-rel}
In the presence of Assumptions~\ref{A2} and~\ref{A3all}.(ii),
the definitions of the scale function $p$ and the speed
measure $m$ imply that
\be
\int _0^x b(s) k(s) \, m(\di s) = \int _0^x k(s) \biggl( \frac{1}{p'}
\biggr) ' (s) \, \di s = \frac{k(x)}{p' (x)} - \frac{1}{2} \int _0^x
\sigma^2 (s) k'(s) \, m(\di s) .
\ee
In turn, these identities and the definition (\ref{K}) of $K$ imply
that
\ben
\int _0^x K(s) \, m(\di s) = \frac{k(x)}{p' (x)} . \label{Kk-rel}
\een
\end{rem}

\begin{rem} \label{rem:lambdas}
In view of Assumption~\ref{A4e}.(ii), we define
\ben
\underline{\lambda} = \lim _{x \uparrow \infty} K(x)
+ h(x) \quad \text{and} \quad \overline{\lambda}
= K(\xi) + h(\xi) , \label{over-under-lam}
\een
and we note that the equation $K(x) + h(x) - \lambda
= 0$ has

{\bf -}
no strictly positive solutions if $\lambda > \overline{\lambda}$,

{\bf -}
two strictly positive solutions if $\lambda \in \mbox{}
\bigl] K(0)+h(0), \, \overline{\lambda} \bigr[$, and

{\bf -}
one strictly positive solution if $\lambda \in \mbox{} \bigl]
\underline{\lambda}, \, K(0)+h(0) \bigr]$ or $\lambda =
\overline{\lambda}$

\noindent
(see also Figure~1).
In particular, there exists a unique function $\varrho $
such that
\ben
\xi < \varrho (\lambda) \quad \text{and} \quad
K \bigl( \varrho (\lambda) \bigr) + h \bigl( \varrho
(\lambda) \bigr) - \lambda = 0 \quad \text{for all }
\lambda \in \mbox{} ]\underline{\lambda},
\overline{\lambda}[ . \label{Kbhl-rho}
\een
Furthermore, this function is such that
\begin{gather}
K(x) + h(x) - \lambda \begin{cases} > 0 , &
\text{for all } x \in \bigl[\xi , \varrho (\lambda)
\bigr[ , \\ < 0 , & \text{for all } x > \varrho (\lambda)
, \end{cases} \label{Kbhl-rho-prop} \\
\text{and} \quad
K' \bigl( \varrho (\lambda) \bigr) + h' \bigl( \varrho
(\lambda) \bigr) < 0 \quad \text{for all } \lambda \in
\mbox{} ]\underline{\lambda}, \overline{\lambda}[ .
\label{Kbhl-der-rho}
\end{gather}
On the other hand, there is a unique function
$\uprho$ such that
\ben
0 < \uprho (\lambda) < \xi \quad \text{and} \quad K
\bigl( \uprho (\lambda) \bigr) + h \bigl( \uprho (\lambda)
\bigr) - \lambda = 0 \quad \text{for all } \lambda \in
\mbox{} \bigl] K(0)+h(0), \, \overline{\lambda} \bigr[ .
\label{Kbhl-uprho}
\een
Given any $\lambda \in \mbox{} \bigl] K(0)+h(0),
\, \overline{\lambda} \bigr[$, this function is such
that
\ben
K(x) + h(x) - \lambda < 0 \quad \text{for all }
x \in \mbox{} \bigl] 0, \uprho (\lambda) \bigr[ .
\label{Kbhl-uprho-prop}
\een
\end{rem}


\begin{figure}[!tbp]
  \centering
    \includegraphics[width=110mm]{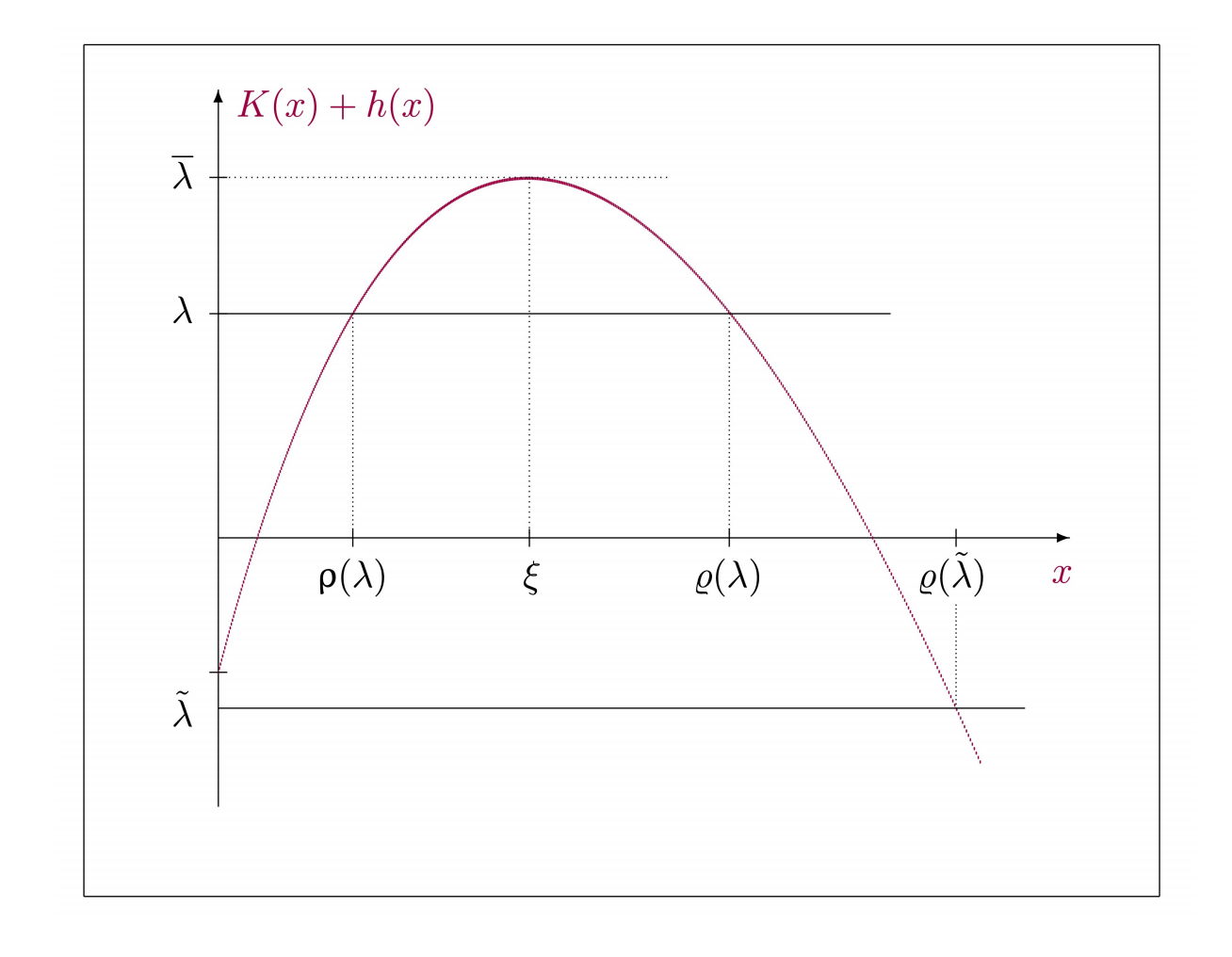}
    \caption{Notation associated with the graph
    of the function $K+h$.}
    \label{fig1}
\end{figure}

We conclude this section with the following three examples.

\begin{ex} \label{ex:VPdif}
Suppose that the uncontrolled dynamics of the state process
are modelled by the SDE
\be
\di X_t = \upkappa (\upgamma - X_t) X_t \, \di t + \upsigma
X_t^\ell \, \di W_t , \quad X_0 = x  > 0 ,
\ee
for some strictly positive constants $\upkappa$, $\upgamma$,
$\upsigma$ and $\ell \in [1, \frac{3}{2}]$.
Note that the celebrated stochastic Verhulst-Pearl logistic
model of population growth arises in the special case
$\ell = 1$.
The derivative of the scale function admits the expression
\begin{gather}
p' (x) = x ^{-2\upkappa \upgamma / \upsigma ^2}
\exp \biggl( \frac{2\upkappa}{\upsigma ^2} (x - 1) \biggr) ,
\quad \text{if } \ell = 1 , \nonumber \\
p' (x) = \exp \biggl( \frac{\upkappa \upgamma}
{(\ell - 1)\upsigma ^2} \bigl( x^{-2(\ell -1)} - 1 \bigr) +
\frac{2\upkappa} {(3-2\ell) \upsigma ^2} \bigl( x^{3-2\ell}
- 1 \bigr) \biggr) , \quad \text{if } \ell \in \mbox{} ]1, 1.5[
, \nonumber \\
\text{and} \quad
p' (x) = x^{2\upkappa / \upsigma ^2}
\exp \biggl( \frac{2\upkappa \upgamma}{\upsigma ^2}
(x^{-1} - 1) \biggr) , \quad \text{if } \ell = 1.5 . \nonumber
\end{gather}
Assumptions~\ref{A1}--\ref{A12+} hold true
if $\ell \in \mbox{} ]1, \frac{3}{2}]$ or if $\ell = 1$
and $\upkappa \upgamma - \frac{1}{2} \sigma^2
> 0$.
Furthermore, if $r$, $k$ are constant and either
(a) $h=0$ and $r < \upkappa \upgamma$,\footnote{The
inequality $r < \upkappa \upgamma$ is needed
for Assumption~\ref{A4d}.(iv) to hold true.
In the absence of this condition, the harvesting
strategy that drives the species to extinction at
time 0 is optimal for the discounted version of the
problem (see Alvarez and
Shepp~\cite[Theorem~1.(i)]{AS}).}
or (b) $h$ is a strictly concave function satisfying the
Inada conditions
\ben
\lim _{x \downarrow 0} h'(x) = \infty \quad \text{and}
\quad \lim _{x \uparrow \infty} h'(x) = 0 , \label{Inada}
\een
as well as the inequality $h(0) > {-\infty}$, then all
of the conditions in Assumptions~\ref{A3all}--\ref{A4e}
are satisfied.
\end{ex}

\begin{ex}
Suppose that the uncontrolled dynamics of the state process
are modelled by the SDE
\be
\di X_t = \biggl( \upkappa \upgamma + \frac{1}{2}
\sigma ^2 - \upkappa \ln (X_t) \biggr) X_t \, \di t +
\upsigma X_t \, \di W_t , \quad X_0 = x  > 0 ,
\ee
for some constants $\upkappa, \upgamma, \upsigma
> 0$, namely, the logarithm of the uncontrolled state
process is the Ornstein-Uhlenbeck process given by
\be
\di \ln (X_t) = \upkappa \bigl( \upgamma - \ln (X_t) \bigr)
\, \di t + \upsigma \, \di W_t , \quad \ln (X_0) = \ln (x)
\in \bbr .
\ee
In this case, the derivative of scale function admits the
expression
\begin{align}
p' (x) & = x^{\frac{\upkappa}{\upsigma^2} \ln (x) -
\frac{2\upkappa \upgamma}{\upsigma^2} - 1} \nonumber
\end{align}
and all of Assumptions~\ref{A1}--\ref{A12+} hold true.
Furthermore, if $r$, $k$ are constant and either $h=0$
or $h$ is a strictly concave function satisfying the
Inada conditions (\ref{Inada}), as well as the
inequality $h(0) > -\infty$, then all of the conditions in
Assumptions~\ref{A3all}--\ref{A4e} are satisfied.
\end{ex}

\begin{ex}
Suppose that the uncontrolled dynamics of the state process
are modelled by the SDE
\be
\di X_t = \upkappa (\upgamma - X_t) \, \di t + \upsigma
X_t^\ell \, \di W_t , \quad X_0 = x  > 0 ,
\ee
for some strictly positive constants $\upkappa$, $\upgamma$,
$\upsigma$ and $\ell \in [ \frac{1}{2} , 1]$.
Note that, in the special case that arises for $\ell = \frac{1}{2}$
and $\upkappa \upgamma - \frac{1}{2} \upsigma^2 > 0$,
the process $X$ identifies with the short rate
process in the Cox-Ingersoll-Ross interest rate model.
The derivative of the scale function admits the expression
\begin{gather}
p' (x) = x^{-2\upkappa \upgamma / \upsigma ^2}
\exp \biggl( \frac{2\upkappa}{\upsigma ^2} (x - 1) \biggr) ,
\quad \text{if } \ell = 0.5 , \nonumber \\
p' (x) = \exp \biggl( \frac{2 \upkappa \upgamma}
{(2\ell - 1)\upsigma ^2} \bigl( x^{-(2\ell -1)} - 1 \bigr) +
\frac{\upkappa} {(1-\ell) \upsigma ^2} \bigl( x^{2(1-\ell)}
- 1 \bigr) \biggr) , \quad \text{if } \ell \in \mbox{} ]
0.5 , 1[ , \nonumber \\
\text{and} \quad
p' (x) = x^{2\upkappa / \upsigma ^2}
\exp \biggl( \frac{2\upkappa \upgamma}{\upsigma ^2}
(x^{-1} - 1) \biggr) , \quad \text{if } \ell = 1 .\nonumber
\end{gather}
Assumptions~\ref{A1}--\ref{A12+} hold true
if $\ell \in \mbox{} ]\frac{1}{2} , 1]$ or if $\ell = \frac{1}{2}$
and $k\upgamma - \frac{1}{2} \sigma^2 > 0$.
Furthermore, if $r$, $k$ are constant and $h$ is a
strictly concave function satisfying the
Inada conditions (\ref{Inada}), as well as the
inequality $h(0) > -\infty$, then all of the conditions
in Assumptions~\ref{A3all}--\ref{A4e} are satisfied.
\end{ex}

\section{The HJB equations of the control problems}
\label{sec:HJB}

We will solve the discounted control problem by deriving
a suitable $C^2$ solution $w$ to the HJB equation
\ben
\max \biggl\{ \frac{1}{2} \sigma ^2 (x) w''(x) + b(x) w'(x)
- r(x) w(x) + h(x) , \ k(x) - w'(x) \biggr\} = 0 \label{HJBd}
\een
that identifies with the control problem's value function.
In particular, we will construct a solution to this HJB
equation such that
\be
\sup_{\zeta \in \acal} I_x (\zeta) = w(x) \quad \text{for all }
x > 0 .
\ee
On the other hand, we will solve both of the ergodic control
problems by constructing a $C^2$ function $w$ and finding
a constant $\lambda$ such that the HJB equation
\ben
\max \biggl\{ \frac{1}{2} \sigma ^2 (x) w''(x) + b(x) w'(x)
+ h(x) - \lambda , \ k(x) - w'(x) \biggr\} = 0 \label{HJBe}
\een
holds true for all $x>0$.
Given such a solution, we will prove that
\be
\sup_{\zeta \in \acal} \JE_x (\zeta) = \lambda \quad
\text{and} \quad \sup_{\zeta \in \acal} \JP_x (\zeta)
= \lambda \quad \text{for all } x > 0 .
\ee

Given suitable solutions to the HJB equations, the optimal
strategies can be characterised as follows.
In the discounted control problem, the controller should wait
and take no action for as long as the state process $X$
takes values in the interior of the set in which the ODE
\ben
\frac{1}{2} \sigma ^2 (x) w''(x) + b(x) w'(x) - r(x) w(x)
+ h(x) = 0 \label{ODE1}
\een
is satisfied.
On the other hand, the ``no-action'' region of the ergodic
control problems is the interior of the set in which the ODE
\ben
\frac{1}{2} \sigma^2 (x) w''(x) + b(x) w'(x) + h(x) - \lambda
= 0 \label{ODE2}
\een
is satisfied.
In all problems, the controller should take the minimal action
required so that the state process is kept outside the interior
of the set defined by $w'(x) = k(x)$ at all times.

We are going to establish that, in all of the problems that
we consider, the optimal strategy takes the following
qualitative form.
There exists a point $\beta$ in the state space
$]0, \infty[$ such that it is optimal to push in an impulsive
way the state process down to level $\beta$ if the initial
state $x$ is strictly greater than $\beta$ and otherwise take
minimal action so that the state process $X$ is kept inside
the set $]0, \beta ]$ at all times, which
amounts at reflecting $X$ in $\beta$ in the negative
direction.
In view of the discussion in the previous paragraph, the
optimality of such a strategy is associated with a solution
$w$ to the HJB equation (\ref{HJBd}) such that
\begin{align}
\frac{1}{2} \sigma^2 (x) w''(x) + b(x) w'(x) - r(x) w(x)
+ h(x) = 0 , & \quad \text{for } x \in \mbox{} ]0, \beta[
, \nonumber \\
w '(x) = k(x) , & \quad \text{for } x \in [\beta, \infty[ ,
\nonumber
\end{align}
and a solution $(w , \lambda)$ to the HJB equation
(\ref{HJBe}) such that
\begin{align}
\frac{1}{2} \sigma^2 (x) w''(x) + b(x) w'(x) + h(x) -
\lambda = 0 , & \quad \text{for } x \in \mbox{} ]0, \beta[
, \nonumber \\
w '(x) = k(x) , & \quad \text{for } x \in [\beta, \infty[ .
\nonumber
\end{align}
In both cases, we will  determine the free-boundary point
$\beta$ using the so-called ``smooth pasting condition"
of singular stochastic control, which requires that $w$
should be $C^2$, in particular, at the free-boundary point
$\beta$.
This condition suggests the free-boundary equations
\ben
\lim_{x \uparrow \beta} w'(x) = k(\beta) \quad \text{and} \quad
\lim_{x \uparrow \beta} w''(x) = k'(\beta) . \label{cond}
\een

\section{The solution to the discounted harvesting problem}
\label{sec:disc-sol}

It is well-known that, in the presence of
Assumptions~\ref{A1}, \ref{A2}, \ref{A3all}.(i)
and~\ref{A4d}.(i)-(ii), every solution to the ODE
(\ref{ODE1}) is given by
\ben
w(x) = A \varphi (x) + B \psi (x) + R_h (x) , \label{w-cont}
\een
for some constants $A, B \in \bbr$.
The functions $\varphi , \psi : \mbox{} ]0, \infty[
\mbox{} \rightarrow \bbr$ are $C^2$, unique up to
a multiplicative constant and satisfying
\begin{gather}
0 < \varphi (x) \quad \text{and} \quad \varphi' (x) < 0 \quad
\text{for all } x>0 , \nonumber \\
0 < \psi (x) \quad \text{and} \quad \psi' (x) > 0 \quad
\text{for all } x>0 \nonumber \\
\text{and} \quad
\lim _{x \downarrow 0} \varphi (x)  = \lim _{x \uparrow \infty}
\psi (x) = \infty , \nonumber
\end{gather}
as well as
\begin{gather}
\varphi (x) \psi' (x) - \varphi' (x) \psi (x) = \bigl(
\varphi (1) \psi' (1) - \varphi' (1) \psi (1) \bigr) p' (x)
=: C p' (x) , \label{C} \\
\frac{\varphi (y)}{\varphi (x)}
= \EXP_y \bigl[ e^{-\Lambda_{T_x}} \bigr]
\quad \text{and} \quad \frac{\psi (x)}{\psi (y)}
= \EXP_x \bigl[ e^{-\Lambda_{T_y}} \bigr]
, \label{phi-psi-defn}
\end{gather}
for all $0 < x \leq y$,
where $\Lambda$ is defined by (\ref{Lam})
for $X^\zeta$ being the solution $X$ to the SDE
(\ref{SDE-0}) and $T_y$ is the first hitting time of
$\{ y \}$ by $X$
(e.g., see Borodin and Salminen~\cite[II.10]{BS}).
The function $R_h: \mbox{} ]0, \infty[ \mbox{} \rightarrow
\bbr$ is $C^2$ and admits the expressions
\begin{align}
\begin{split}
R_h (x) & = \EXP_x \biggl[ \int _0^\infty e^{-\Lambda_t}
h(X_t) \, \di t \biggr] \\
& = \varphi (x) \int _0^x h(s) \Psi (s) \, \di s + \psi (x)
\int _x^\infty h(s) \Phi (s) \, \di s ,
\end{split}
\label{Rh}
\end{align}
in which, $X$ is the solution to the SDE (\ref{SDE-0}),
\begin{align}
\Phi (x) & = \frac{2 \varphi (x)}{C \sigma^2 (x) p'(x)}
= \frac{1}{C} \varphi (x) \frac{m (\di x)}{\di x}
\label{Phi} \\
\text{and} \quad
\Psi (x) & = \frac{2 \psi (x)}{C \sigma^2 (x) p'(x)}
= \frac{1}{C} \psi (x) \frac{m (\di x)}{\di x} .
\label{Psi}
\end{align}
In view of the definition (\ref{Theta}) of $\Theta$,
the conditions in Assumptions~\ref{A4d}.(iii)
imply that
\ben
R_h (x) = R_{\Theta + h} (x) + \int _0^x k(s) \, \di s
- \Theta_\infty \psi (x) , \label{Rh-RTheta}
\een
where $R_{\Theta + h}$ is given by (\ref{Rh}) with
$\Theta + h$ in place of $h$ and
\ben
\Theta _\infty = \lim _{x \uparrow \infty} \frac{1}{\psi (x)}
\int _0^x k(s) \, \di s \in \bbr _+ . \label{Thetaoo}
\een
All of these claims follow from the results in Lamberton
and Zervos~\cite[Section 4]{LZ} (see also Liu and
Zervos~\cite[Section 3]{LZ24}).

In (\ref{w-cont}), we choose $A=0$ because, otherwise,
a solution to the HJB equation (\ref{HJBd}) of the
form we have discussed in the previous section
cannot satisfy the so-called ``transversality condition'',
which is required for a solution $w$ to the HJB equation
to identify with the control problem's value function.
In view of this observation and the identity
(\ref{Rh-RTheta}), we consider the solution to the ODE
(\ref{ODE1}) that is given by
\be
w(x) = R_{\Theta + h} (x) + \int _0^x k(s) \, \di s +
(B - \Theta _\infty) \psi (x) ,
\ee
for some constant $B \in \bbr$.
Substituting $B$ for its unique choice for which this
function satisfies the first boundary condition in
(\ref{cond}), we obtain
\begin{align}
\begin{split}
w(x) & = \int _0^x k(s) \, \di s + R_{\Theta + h} (x) -
\frac{R_{\Theta + h}' (\beta)}{\psi' (\beta)} \psi (x) \\
& = R_h (x) + \biggl( \Theta _\infty -
\frac{R_{\Theta + h}' (\beta)}{\psi' (\beta)} \biggr)
\psi (x), \quad \text{for } x \in \mbox{} ]0, \beta] .
\end{split}
\label{OS3}
\end{align}
Furthermore, this function satisfies the second
boundary condition in (\ref{cond}) if and only if
$\beta > 0$ is a solution to the algebraic equation
\ben
R_{\Theta + h}'' (\beta) -
\frac{R_{\Theta + h}' (\beta)}{\psi' (\beta)} \psi'' (\beta)
= 0 . \label{F}
\een

\begin{lemma} \label{lem-existGamma}
Consider the discounted control problem formulated
in Section~\ref{sec:pr-form}.
There exists a unique point $\beta^\star \in \mbox{}
]\xi, \infty[$ satisfying equation (\ref{F}), which is
equivalent to
\ben
\int _0^\beta \bigl( \Theta (s) + h(s) \bigr) \Psi (s)
\, \di s = \bigl( \Theta (\beta) + h(\beta) \bigr) \frac{1}
{r(\beta)} \int _0^\beta r(s) \Psi (s) \, \di s . \label{FA}
\een
Furthermore, the function $w$ defined by
\ben
w(x) = \begin{cases} R_h (x) + \Bigl( \Theta_\infty -
\frac{R_{\Theta + h}' (\beta^\star)}{\psi' (\beta^\star)}
\Bigr) \psi(x) , & \text{for } x \in \mbox{} ]0, \beta^\star]
, \\ w(\beta^\star) + \int _{\beta^\star}^x k(s)
\, \di s , & \text{for } x \in \mbox{} ]\beta^\star, \infty[ ,
\end{cases} \label{w-sol}
\een
is a $C^2$ solution to the HJB equation (\ref{HJBd})
that is bounded from below.
\end{lemma}

\begin{proof}
In the presence of the assumptions that we have made,
Lemma~4 in Liu and Zervos~\cite{LZ24} proves that
that there exists a unique point $\beta^\star \in \mbox{}
]\xi, \infty[$ (denoted by $\underline{x}$ in that lemma)
such that
\ben
\frac{\di}{\di x} \frac{R'_{\Theta + h} (x)}{\psi' (x)} \begin{cases}
< 0 & \text{for all } x \in \mbox{} ]0, \beta^\star[ , \\ > 0 &
\text{for all } x \in \mbox{} ]\beta^\star, \infty[ , \end{cases}
\quad \text{and} \quad
\lim _{x \uparrow \infty} \frac{R'_{\Theta + h} (x)}{\psi' (x)}
= 0 . \label{beta*}
\een
The first set of these inequalities and the continuity of
the derivative $(R'_{\Theta + h} / \psi')'$ imply that
$\beta^\star$ is the unique solution to equation
(\ref{F}).

To establish the equivalence of (\ref{F}) and (\ref{FA}),
we use the definitions (\ref{Rh})--(\ref{Psi}) to
see that (\ref{F}) is equivalent to
\be
\int _0^\beta \bigl( \Theta (s) + h(s) \bigr) \Psi (s) \,
\di s = \bigl( \Theta (\beta) + h(\beta) \bigr) \frac{\psi' (\beta)}
{C r(\beta) p'  (\beta)} .
\ee
In the presence of Assumption~\ref{A2}, 0 is an
inaccessible boundary point of the solution $X$ to the
SDE (\ref{SDE-0}), therefore, $\lim _{x \downarrow 0}
\psi' (x) / p' (x) = 0$ (e.g., see Borodin and
Salminen~\cite[Section~II.10]{BS}).
Combining this limit with the calculation
\be
\frac{\di}{\di x} \frac{\psi' (x)}{p' (x)} = \frac{2}{\sigma^2 (x) p'(x)}
\biggl( \frac{1}{2} \sigma^2 (x) \psi'' (x) + b(x) \psi' (x) \biggr)
= \frac{2  r(x) \psi (x)}{\sigma^2 (x) p'(x)} = C r(x) \Psi (x) ,
\ee
we obtain
\ben
\int _0^x r(s) \Psi (s) \, \di s = \frac{\psi' (x)}{C p' (x)}
\quad \text{for all } x>0 , \label{psi/p'-int}
\een
and the equivalence of (\ref{F}) and (\ref{FA})
follows.

The assumption that $h$ is bounded from below implies
that
\ben
\inf_{x>0} R_h (x) = \inf _{x>0} \EXP_x \biggl[ \int
_0^\infty e^{-\Lambda_t} h(X_t) \, \di t \biggr] \geq
\inf _{x>0} h(x) \EXP_x \biggl[ \int _0^\infty
e^{-\Lambda_t} \, \di t \biggr] > - \infty . \label{w-BFB}
\een
Combining this observation with the fact that
$\Theta _\infty \in \bbr_+$ (see the remark after
(\ref{Thetaoo})) and the positivity of $k$, we can see
that $w$ is bounded from below.

By construction, we will establish all of the lemma's
other claims if we prove that
\begin{align}
k(x) - w'(x) \leq 0 & \quad \text{for all } x \in \mbox{}
]0, \beta^\star[ \label{HJBd-ineq1} \\
\text{and} \quad \frac{1}{2} \sigma^2 (x) w''(x) + b(x)
w'(x) - r(x) w(x) + h(x) \leq 0 & \quad \text{for all }
x \in \mbox{} ]\beta^\star, \infty[ . \label{HJBd-ineq2}
\end{align}
The inequality (\ref{HJBd-ineq1}) holds true because
the first expression of $w$ in (\ref{OS3}) and (\ref{beta*})
imply that
\be
k(x) - w'(x) = \psi' (x) \biggl(
\frac{R_{\Theta + h}' (\beta^\star)}{\psi' (\beta^\star)}
- \frac{R_{\Theta + h}' (x)}{\psi' (x)} \biggr)
< 0 \quad \text{for all } x \in \mbox{} ]0, \beta^\star[ .
\ee
To show (\ref{HJBd-ineq2}), we first use the expression
of $w$ given by (\ref{w-sol}) for $x > \beta^\star$,
the definition the definition (\ref{Theta}) of $\Theta$
and the first expression in (\ref{OS3}) to calculate
\begin{align}
\frac{1}{2} \sigma^2 (x) w''(x) + b(x) & w'(x) - r(x) w(x)
+ h(x) \nonumber \\
& = \Theta (x) + h(x) - r(x) \biggl( w(\beta^\star) -
\int _0^{\beta^\star} k(s) \, \di s \biggr)
\label{HJBd-ineq2+} \\
& = \Theta (x) + h(x) -r(x) G_{\Theta + h} (\beta^\star)
, \nonumber
\end{align}
where
\begin{align}
G_{\Theta + h} (x) & = R_{\Theta + h} (x) -
\frac{R_{\Theta + h}' (x)}{\psi' (x)} \psi (x) =
\frac{C p' (x)}{\psi' (x)} \int _0^x \bigl( \Theta (s)
+ h(s) \bigr) \Psi (s) \, \di s . \label{G-Theta}
\end{align}
In view of the calculations
\begin{align}
G_{\Theta + h}' (x) & = - \psi (x) \frac{\di}{\di x}
\frac{R_{\Theta + h}' (x)}{\psi' (x)} \nonumber \\
& = - \frac{2C r(x) p'(x) \psi (x)}{( \sigma (x) \psi' (x) )^2}
\biggl( \int _0^x \bigl( \Theta (s) + h(s) \bigr)
\Psi (s) \, \di s - \frac{\Theta (x) + h (x)}{r(x)}
\frac{\psi' (x)}{C p'(x)} \biggr) , \nonumber
\end{align}
and the inequalities (\ref{beta*}), we can see that
\begin{gather}
G_{\Theta + h} (x) < G_{\Theta + h} (\beta^\star)
\label{HJB-ineq221} \\
\text{and} \quad \int _0^x \bigl( \Theta (s) + h(s)
\bigr) \Psi (s) \, \di s > \frac{\Theta (x) + h(x)}{r(x)}
\frac{\psi' (x)}{C p'(x)} \label{HJB-ineq222}
\end{gather}
for all $x > \beta^\star$.
The inequality (\ref{HJB-ineq222}) and the second
expression of $G_{\Theta+h}$ in (\ref{G-Theta})
imply that
\be
G_{\Theta + h} (x) = \frac{C p'(x)}{\psi' (x)} \int _0^x
\bigl( \Theta (s) + h(s) \bigr) \Psi (s) \, \di s >
\frac{\Theta (x) + h(x)}{r(x)} \quad \text{for all }
x > \beta^\star .
\ee
However, this result, (\ref{HJBd-ineq2+}) and the
inequality (\ref{HJB-ineq221}) yield
\be
\frac{1}{2} \sigma^2 (x) w''(x) + b(x) w'(x) - r(x) w(x)
+ h(x) < r(x) \biggl( \frac{\Theta (x) + h(x)}{r(x)} -
G_{\Theta + h} (x) \biggr) < 0
\ee
for all $x > \beta^\star$, and (\ref{HJBd-ineq2})
follows.
\end{proof}

\begin{theorem} \label{thm:SolPb-lem}
Consider the discounted control problem formulated in
Section~\ref{sec:pr-form}.
If the point $\beta^\star \in \mbox{} ]\xi, \infty[$ and
the function $w$ are as in Lemma~\ref{lem-existGamma},
then
\ben
w(x) = \sup_{\zeta \in \acal} J_x(\zeta) \quad
\text{for all } x > 0 , \label{supI=w}
\een
while the harvesting strategy $\zeta^\star \in \acal$
that has a jump of size $\Delta \zeta _0^\star
= (x-\beta^\star)^+$ at time 0 and then reflects the
state process $X^\star$ at the level $\beta^\star$
in the negative direction is optimal.
\end{theorem}

\begin{proof}
Fix any initial value $x>0$, consider any admissible controlled
process $\zeta \in \acal$ and denote by $X^\zeta$ the
associated solution to the SDE (\ref{SDE}).
Using It\^{o}'s formula,  we obtain
\begin{align}
e^{-\Lambda _T^\zeta} & w(X _T^\zeta) \nonumber \\
= \mbox{} & w(x) +
\int _0^T e^{-\Lambda _t^\zeta} \biggl( \frac{1}{2} \sigma^2
(X_t^\zeta) w'' (X_t^\zeta) + b(X_t^\zeta) w' (X_t^\zeta)
- r(X_t^\zeta) w(X_t^\zeta) \biggr) \, \di t + M_T^\zeta
\nonumber \\
& - \int _{[0,T]} e^{-\Lambda _t^\zeta} w' (X_{t-}^\zeta)
\, \di \zeta_t + \sum _{0 \leq t \leq T} e^{-\Lambda _t^\zeta}
\Bigl( w(X_t^\zeta) - w(X_{t-}^\zeta) - w' (X_{t-}^\zeta)
\, \Delta X_t^\zeta \Bigr) , \nonumber
\end{align}
where
\be
M_T^\zeta = \int _0^T e^{-\Lambda _t^\zeta} \sigma
(X_t^\zeta) w' (X_t^\zeta) \, \di W_t .
\ee
Since $\Delta X_t^\zeta = X_t^\zeta - X^\zeta_{t-}
= - \Delta \zeta_t \leq 0$, we can see that
\begin{align}
w (X_t^\zeta) - w (X_{t-}^\zeta) + \int _0^{\Delta \zeta_t}
& k(X_{t-}^\zeta - u) \, \di u = \int _{X_t^\zeta}^{X_{t-}^\zeta}
\bigl( k(u) - w' (u) \bigr) \, \di u  . \nonumber
\end{align}
In view of the facts that $\zeta^c$ is an increasing
process, $X_t^\zeta < X_{t-}^\zeta$ and $w$
satisfies the HJB equation (\ref{HJBd}), we can
see that these observations imply that
\begin{align}
\int _0^T e^{-\Lambda _t^\zeta} & h(X_t^\zeta) \, \di t
+ \int _0^T e^{-\Lambda_t^\zeta} k(X_t^\zeta) \circ
\di \zeta_t \nonumber \\
= \mbox{} & w(x) - e^{-\Lambda _T^\zeta}
w(X_T^\zeta) + M_T^\zeta \nonumber \\
& + \int _0^T e^{-\Lambda _t^\zeta} \biggl( \frac{1}{2}
\sigma ^2 (X_t^\zeta) w'' (X_t^\zeta) + b(X_t^\zeta)
w' (X_t^\zeta) - r(X_t^\zeta) w(X_t^\zeta) + h(X_t^\zeta)
\biggr) \, \di t \label{VTdineq1} \\
& + \int _0^T e^{-\Lambda _t^\zeta} \bigl( k(X_t^\zeta)
- w' (X_t^\zeta) \bigr) \, \di \zeta_t^\cc
+ \sum _{0 \leq t \leq T} e^{-\Lambda _t^\zeta}
\int _{X_t^\zeta}^{X_{t-}^\zeta} \bigl( k(u) - w' (u) \bigr)
\, \di u \nonumber \\
\leq \mbox{} & w(x) - e^{-\Lambda _T^\zeta}
w(X_T^\zeta) + M_T^\zeta . \nonumber
\end{align}

We next consider any sequence $(\tau_n)$ of bounded
localising stopping times for the local martingale
$M^\zeta$.
Recalling that $h$ and $w$ are both bounded from below,
$k$ is positive and $\zeta$ is an increasing process,
we use the dominated and the monotone convergence
theorems to observe that (\ref{VTdineq1}) implies that
\begin{align}
\begin{split}
I_x (\zeta) & = \lim _{n \uparrow \infty} \EXP_x \Biggl[
\int _0^{\tau_n} e^{-\Lambda _t^\zeta} h (X_t^\zeta) \,
\di t + \int _0^{\tau_n} e^{-\Lambda_t^\zeta} k(X_t^\zeta)
\circ \di \zeta_t \biggr] \\
& \leq \lim _{n \uparrow \infty} \EXP_x \Bigl[ w(x) +
e^{-\Lambda _{\tau_n}^\zeta} w^- (X_{\tau_n}^\zeta)
\Bigr] = w(x) ,
\end{split}
\label{I<w}
\end{align}
where $w^- (x) = - \min \bigl\{ 0, w(x) \bigr\}$.

Consider the harvesting strategy $\zeta^\star \in \acal$
that is as in the statement of the theorem: such a strategy
indeed exists (see Tanaka~\cite[Theorem~4.1]{Tan}).
This strategy is such that (\ref{VTdineq1}) holds true
with equality, namely,
\ben
\int _0^T e^{-\Lambda _t^\star} h(X_t^\star) \, \di t
+ \int _0^T e^{-\Lambda_t^\star} k(X_t^\star) \circ
\di \zeta_t^\star = w(x) - e^{-\Lambda _T^\star}
w(X_T^\star) + M_T^\star . \nonumber
\een
Furthermore, the processes $h(X^\star)$,
$k(X^\star)$ and $w(X^\star)$ are all bounded
because $X_t^\star \in \mbox{} ]0, \beta^\star]$
for all $t>0$.
In view of these observations, we can use the
dominated convergence theorem to obtain
\begin{align}
I_x (\zeta^\star) & = \lim _{n \uparrow \infty} \EXP_x
\Biggl[ \int _0^{\tau_n^\star} e^{-\Lambda _t^\star}
h (X_t^\star) \, \di t + \int _0^{\tau_n^\star}
e^{-\Lambda_t^\star} k(X_t^\star) \circ \di \zeta_t^\star
\biggr] \nonumber \\
& = \lim _{n \uparrow \infty} \EXP_x \Bigl[ w(x) -
e^{-\Lambda _{\tau_n^\star}^\star}
w \bigl( X_{\tau_n^\star}^\star \bigr) \Bigr] = w(x)
, \nonumber
\end{align}
where $(\tau_n^\star)$ is a sequence of bounded
localising stopping times for the local martingale
$M^\star$.
However, these identities and (\ref{I<w}) imply
(\ref{supI=w}) as well as the optimality of
$\zeta^\star$.
\end{proof}

\section{The solution to the ergodic problem's HJB equation}
\label{sec:HJB-sol}

In view of (\ref{Kk-rel}) in Remark~\ref{rem:Kk-rel},
we can verify that a solution to the ODE (\ref{ODE2}) is given
by
\ben
w'(x) = w'(x; \lambda) = p'(x) \int _0^x \bigl( \lambda
- h(s) \bigr) \, m(\di s) = k(x) - p' (x) \Xi (x, \lambda) ,
\label{w'1}
\een
where
\ben
\Xi (x, \lambda) = \int _0^x \bigl( K(s) + h(s) - \lambda
\bigr) \, m(\di s) . \label{Xi}
\een
This function satisfies the boundary conditions
(\ref{cond}) if and only if
\ben
\Xi (\beta, \lambda) = 0 \quad \text{and} \quad
K (\beta) + h (\beta) - \lambda = 0 . \label{lam-beta-range}
\een

In the next result, we derive a unique solution to the
system of equations in (\ref{lam-beta-range}) as well
as a solution to the HJB equation (\ref{HJBe}).
It turns out that the solution to the HJB equation
(\ref{HJBe}) may be unbounded from below (see
Remark~\ref{rem:why(III)} at the end of the section),
which gives rise to a non-trivial complication in the
verification arguments we use for the proof of
Theorem~\ref{VT-th}.
The introduction of the auxiliary function $w_\lambda$
in part~(III) of the next result provides a way to
overcome this problem. \blue
To overcome the same complication,
Cohen, Hening and Sun~\cite{CHS22} followed
a similar but much more complex approach
in the presence of stronger assumptions on
the problem data.
In particular, the auxiliary function that they construct
for the proof of their Proposition~3.4 is $C^2$.
In contrast, the function $w_\lambda$
that we consider is $C^1$ but loses
the $C^2$ regularity at two points.
\black

\begin{proposition} \label{SolPb-lem}
In the presence of Assumptions~\ref{A1}, \ref{A2},
\ref{A3all} and~\ref{A4e}, the following statements
hold true:
\smallskip

\noindent {\rm (I)}
There exists a unique pair $(\beta^\star, \lambda^\star)$
with $\beta^\star > 0$ satisfying the system of equations
(\ref{lam-beta-range}).
This pair is such that
\ben
K(0) + h(0) < \lambda^\star = \frac{1}
{m \bigl( ]0, \beta^\star[ \bigr)}
\int _0^{\beta^\star} \bigl( K(s) + h(s) \bigr) \, m(\di s)
< \overline{\lambda}
\text{ and }
\beta^\star = \varrho (\lambda^\star) ,
\label{lambda*location}
\een
where $\overline{\lambda}$ and $\varrho$ are given by
(\ref{over-under-lam}) and (\ref{Kbhl-rho}).
\smallskip

\noindent {\rm (II)}
The unique, modulo an additive constant, function $w$
that is defined by
\ben
w' (x) = \begin{cases} k(x) - p' (x) \Xi (x, \lambda^\star)
, & \text{for } x \in \mbox{} ]0, \beta^\star[ , \\ k(x) , &
\text{for } x \geq \beta^\star , \end{cases} \label{w*-prop}
\een
is a $C^2$ solution to the HJB equation (\ref{HJBe}).
\smallskip

\noindent {\rm (III)}
Given any $\lambda \in \mbox{} ]\lambda^\star ,
\overline{\lambda}[$, there exists a point $\alpha
(\lambda) \in \mbox{} \bigl] 0 , \uprho (\lambda) \bigr[$,
where $\uprho$ is as in (\ref{Kbhl-uprho}), such
that the unique, modulo an additive constant,
function $w_\lambda$ that is defined by
\begin{gather}
w_\lambda' (x) = \int _{\alpha (\lambda)}^x \bigl(
K(s) + h(s) - \lambda \bigr) \, m(\di s) , \quad \text{for }
x \in \mbox{} \bigl] \alpha (\lambda), \varrho (\lambda)
\bigr[ , \label{w-lam-prop1} \\
\text{and} \quad
w_\lambda' (x) = k(x) , \quad \text{for } x \in
\mbox{} \bigl] 0, \alpha (\lambda) \bigr] \cup
\bigl[ \varrho (\lambda) , \infty \bigr[ ,
\label{w-lam-prop2}
\end{gather}
is $C^1$ in $\bbr_+$ and $C^2$ in $\bbr_+ \setminus
\bigl\{ \alpha (\lambda) , \varrho (\lambda) \bigr\}$,
\ben
\lim _{\lambda \downarrow \lambda^\star}
\alpha (\lambda) = 0 , \quad
\lim _{\lambda \downarrow \lambda^\star}
\varrho (\lambda) = \beta^\star
\quad \text{and} \quad
\lim _{\lambda \downarrow \lambda^\star}
w_\lambda' (x) = w' (x) \text{ for all } x > 0 .
\label{wlam->w}
\een
Furthermore, this function is such that
\begin{gather}
w_\lambda' (x) \geq k(x) , \quad \text{if }
x \in \mbox{} \bigl] \alpha (\lambda), \varrho
(\lambda) \bigr[ , \label{w-lam-ineq1} \\
\frac{1}{2} \sigma^2 (x) w_\lambda''(x) + b(x)
w_\lambda'(x) + h(x) - \lambda \begin{cases}
= 0 , & \text{if } x \in \mbox{} \bigl] \alpha
(\lambda), \varrho (\lambda) \bigr[ , \\
< 0 , & \text{if } x \in \mbox{} \bigl] 0, \alpha
(\lambda) \bigr[ \mbox{} \cup  \mbox{} \bigl]
\varrho (\lambda) , \infty \bigr[ , \end{cases}
\label{w-lam-ineq2} \\
\text{and} \quad
\bigl| w_\lambda' (x) \bigr| \leq C_2 \quad
\text{for all } x > 0 , \label{w'-bound}
\end{gather}
for some constant $C_2 = C_2 (\lambda) > 0$.
\end{proposition}

\begin{proof}
We develop the proof in four steps.
\smallskip

\underline{\em Preliminary results\/}.
Given any $\beta > 0$, the definition (\ref{Xi}) of $\Xi$
implies that
\be
\Lambda (\beta) = \frac{1}{m \bigl( ]0, \beta[ \bigr)}
\int _0^\beta \bigl( K(s) + h(s) \bigr) \, m(\di s)
\ee
is the unique solution to the equation $\Xi (\beta,
\lambda) = 0$.
In light of Assumption~\ref{A4e}.(ii) (see also Figure~1),
a straightforward inspection of the definition of $\Xi$
reveals that this solution is such that one of the following
two cases holds true:
\begin{align}
\begin{split}
\text{(i) } & K(0) + h(0) < \Lambda (\beta) <
\overline{\lambda} \\
\text{or} \quad
\text{(ii) } & \underline{\lambda} < \Lambda (\beta)
\leq K(0) + h(0) \text{ and } \varrho \bigl( \Lambda
(\beta) \bigr) < \beta ,
\end{split}
\label{Lambda(beta)-cases}
\end{align}
where $\underline{\lambda} < \overline{\lambda}$ are
defined by (\ref{over-under-lam}) and $\varrho$ is
introduced by (\ref{Kbhl-rho}).
In particular, we note that $\Lambda (\beta) \in \mbox{}
]\underline{\lambda} , \overline{\lambda} [$, which is the
domain of the function $\varrho$.

Differentiating the identity
\ben
\Xi \bigl( \beta, \Lambda (\beta) \bigr) = 0 , \quad
\text{for } \beta > 0 , \label{Lambda}
\een
which defines $\Lambda$, with respect to $\beta$,
we calculate
\be
\Lambda' (\beta) = \frac{2}
{\sigma^2 (\beta) p' (\beta) m \bigl( ]0 ,\beta[ \bigr)}
\bigl( K(\beta) + h(\beta) - \Lambda (\beta) \bigr) .
\ee
On the other hand, differentiating the identity
\be
K \bigl( \varrho \bigl( \Lambda (\beta) \bigr) \bigr) +
h \bigl( \varrho \bigl( \Lambda (\beta) \bigr) \bigr) -
\Lambda (\beta) = 0 ,
\ee
which follows from (\ref{Kbhl-rho}), with respect to $\beta$,
we derive the expression
\be
\Lambda' (\beta) = \Bigl( K' \bigl( \varrho \bigl( \Lambda
(\beta) \bigr) \bigr) + h' \bigl( \varrho \bigl( \Lambda (\beta)
\bigr) \bigr) \Bigr) \, \frac{\di}{\di \beta} \varrho \bigl(
\Lambda (\beta) \bigr) .
\ee
Combining these calculations, we obtain
\be
\frac{\di}{\di \beta} \varrho \bigl( \Lambda (\beta) \bigr) =
\frac{2 \bigl( K(\beta) + h(\beta) - \Lambda (\beta) \bigr)}
{\sigma^2 (\beta) p' (\beta) m \bigl( ]0, \beta[ \bigl) \bigl(
K' \bigl( \varrho (\Lambda (\beta)) \bigr) + h' \bigl( \varrho
(\Lambda (\beta)) \bigr) \bigr)} .
\ee
In view of this result and the inequality
\be
K' \bigl( \varrho \bigl( \Lambda (\beta) \bigr) \bigr)
+ h' \bigl( \varrho \bigl( \Lambda (\beta) \bigr) \bigr)
< 0 \quad \text{for all } \beta > 0 ,
\ee
which follows from (\ref{Kbhl-der-rho}), we can
see that
\ben
\sgn \biggl( \frac{\di}{\di \beta} \varrho \bigl( \Lambda (\beta)
\bigr) \biggr) = - \sgn \Bigl( K(\beta) + h(\beta) - \Lambda
(\beta) \Bigr) \quad \text{for all } \beta > 0 ,
\label{sgn(rho(L))'}
\een
where $\sgn$ is the sign function defined by
\be
\sgn (x) = \begin{cases} \frac{x}{|x|} , & \text{for } x \neq 0, \\
0 , & \text{for } x = 0 . \end{cases}
\ee

\underline{{\em Proof of} (I).}
In view of (\ref{Lambda}), we can see that there exists a
pair $(\beta^\star, \lambda^\star)$ with $\beta^\star > 0$
satisfying the system of equations (\ref{lam-beta-range})
if and
only if
\ben
K(\beta^\star) + h (\beta^\star) - \lambda^\star = 0
\quad \text{and} \quad
\lambda^\star = \Lambda (\beta^\star) .
\label{lam-bet-eqn0}
\een
The structure of the function $K+h$, which we
have discussed in Remark~\ref{rem:lambdas}
(see also Figure~1), implies that there exists no
$\beta^\star$ satisfying (\ref{lam-bet-eqn0}) if
$\Lambda (\beta^\star)$ is as in case (ii) of 
(\ref{Lambda(beta)-cases}).
We therefore need to show that there exists $\beta^\star
> 0$ such that, if we define $\lambda^\star = \Lambda
(\beta^\star)$, then $\lambda^\star$ satisfies the inequalities (\ref{lambda*location}), and
\ben
\text{either} \quad \beta^\star = \uprho \bigl( \Lambda
(\beta^\star) \bigr) = \uprho (\lambda^\star)
\quad \text{or} \quad
\beta^\star = \varrho \bigl( \Lambda (\beta^\star) \bigr)
= \varrho (\lambda^\star) , \label{lam-bet-eqn1}
\een
where $\varrho$, $\uprho$ are as in (\ref{Kbhl-rho}),
(\ref{Kbhl-uprho}).
Furthermore, the resulting solution $(\beta^\star,
\lambda^\star)$ to the system of equations
(\ref{lam-beta-range}) is unique if and only if only one
of the two equations in (\ref{lam-bet-eqn1}) has a
unique solution and the other one has no solution.

If the equation $\beta = \uprho \bigl( \Lambda (\beta)
\bigr)$ had a solution $\beta^\star > 0$, then
(\ref{Kbhl-uprho-prop}) would imply that
\be
K(s) + h(s) - \Lambda(\beta^\star) < 0 \quad \text{for all }
s < \uprho \bigl( \Lambda (\beta^\star) \bigr) = \beta^\star ,
\ee
which would contradict the identity
\be
\Xi \bigl( \beta^\star, \Lambda (\beta^\star) \bigr)
= \int _0^{\beta^\star} \bigl( K(s) + h(s) - \Lambda
(\beta^\star) \bigr) \, m(\di s) = 0 .
\ee
To establish part~(I) of the theorem, we therefore have to
prove that there exists a unique point $\beta^\star > 0$
such that
\ben
\Lambda (\beta^\star) \in \mbox{} \bigl] K(0) + h(0) ,
\, \overline{\lambda} \bigr[ \quad \text{and} \quad
\beta^\star = \varrho \bigl( \Lambda (\beta^\star) \bigr) .
\label{bet-Lam-*}
\een

To prove that there exists a unique $\beta^\star > 0$
satisfying (\ref{bet-Lam-*}), we first observe that the
inequality in (\ref{Kbhl-rho}) implies that
\ben
\beta < \varrho \bigl( \Lambda (\beta) \bigr) \quad
\text{for all } \beta \leq \xi . \label{beta<rho(beta)}
\een
We next argue by contradiction and we assume that
there is no $\beta^\star > 0$ satisfying the equation in
(\ref{bet-Lam-*}).
In view of (\ref{beta<rho(beta)}) and the continuity
of the functions $\varrho$ and $\Lambda$, we can
see that such an assumption implies that
\ben
\beta < \varrho \bigl( \Lambda (\beta) \bigr) \quad
\text{for all } \beta > \xi . \label{beta<rho(beta)1}
\een
In turn, this inequality and (\ref{Kbhl-rho-prop}) imply
that
\be
K(\beta) + h(\beta) - \Lambda (\beta) > 0 \quad
\text{for all } \beta  > \xi .
\ee
Combining this observation with (\ref{sgn(rho(L))'}), we obtain
$\frac{\di}{\di \beta} \varrho \bigl( \Lambda (\beta) \bigr)
< 0$ for all $\beta > \xi$.
Therefore,
\be
\frac{\di}{\di \beta} \Bigl( \beta - \varrho \bigl( \Lambda
(\beta) \bigr) \Bigr) > 1 \quad \text{for all } \beta > \xi ,
\ee
which contradicts (\ref{beta<rho(beta)1}).
It follows that there exists $\beta^\star > 0$ satisfying
the equation in (\ref{bet-Lam-*}).

To see that the solution $\beta^\star > \xi$ to the equation
in (\ref{bet-Lam-*}) is indeed unique, we note that
(\ref{Kbhl-rho}) implies that $K (\beta) + h(\beta)
- \Lambda (\beta) = 0$ for all $\beta > \xi$ such that
$\beta = \varrho \bigl( \Lambda (\beta) \bigr)$.
This observation and (\ref{sgn(rho(L))'}) imply
\be
\frac{\di}{\di \beta} \bigl( \beta - \varrho (\beta) \bigr)
= 1 \quad \text{for all } \beta > \xi \text{ such that }
\beta = \varrho \bigl( \Lambda (\beta) \bigr).
\ee
Based on this result, we can develop a simple contradiction
argument to show that the equation in (\ref{bet-Lam-*}) has at
most one solution $\beta^\star > \xi$.

We conclude this part of the proof by noting that the
first statement in (\ref{bet-Lam-*}) can be seen by a
straightforward inspection of the equation (\ref{Xi})
that $\bigl( \beta^\star , \Lambda (\beta^\star) \bigr)$
satisfies in light of the identity in (\ref{bet-Lam-*})
and Figure~1.

\underline{{\em Proof of} (II).}
By construction, we will show that the function $w$
given by (\ref{w*-prop}) is a $C^2$ solution to the HJB
equation (\ref{HJBe}) if we prove that
\begin{align}
w '(x) \geq k(x) & \quad \text{for all } x \in \mbox{}
]0, \beta^\star[ \nonumber \\
\text{and} \quad
\frac{1}{2} \sigma^2 (x) w''(x) + b(x) w'(x) + h(x) -
\lambda^\star \leq 0 & \quad \text{for all } x \in \mbox{}
]\beta^\star, \infty[ . \nonumber
\end{align}
In view of the identities $\lambda^\star = \Lambda
(\beta^\star)$ and $\beta^\star = \varrho \bigl( \Lambda
(\beta^\star) \bigr)$, the second of these inequalities
is equivalent to
\be
K(x) + h(x) - \Lambda (\beta^\star) \leq 0 \quad
\text{for all } x > \beta^\star = \varrho \bigl( \Lambda
(\beta^\star) \bigr) ,
\ee
which is true thanks to (\ref{Kbhl-rho-prop}).
On the other hand, the first of these inequalities
follows immediately from the expression of $w'$ in
(\ref{w*-prop}) and the inequalities
\be
\frac{\di}{\di x} \int _0^x \bigl( K(s) +
h(s) - \lambda^\star \bigr) \, m(\di s)
\begin{cases} < 0 & \text{for all } x \in \mbox{} \bigl]
0, \uprho (\lambda^\star) \bigr[ , \\ > 0 & \text{for all }
x \in \mbox{} \bigl] \uprho (\lambda^\star) , \beta^\star
\bigr[ , \end{cases}
\ee
which hold true thanks to the identities $\lambda^\star
= \Lambda (\beta^\star)$ and $\beta^\star = \varrho
\bigl( \Lambda (\beta^\star) \bigr)$, the inequalities
in (\ref{lambda*location}) and Assumption~\ref{A4e}.(ii)
(see also Figure~1).
\smallskip

\underline{{\em Proof of} (III).}
Fix any $\lambda \in \mbox{} ]\lambda^\star ,
\overline{\lambda}[$.
In view of Assumption~\ref{A4e}.(ii) and the properties
of the functions $\varrho$, $\uprho$ in (\ref{Kbhl-rho}),
(\ref{Kbhl-uprho}), we can see that
\begin{align}
\int _0^{\varrho (\lambda)} \bigl( K(s) + h(s) - \lambda
\bigr) \, m(\di s) & < \int _0^{\varrho (\lambda)} \bigl(
K(s) + h(s) - \lambda^\star \bigr) \, m(\di s)
\nonumber \\
& < \int _0^{\beta^\star} \bigl( K(s) + h(s) - \lambda^\star
\bigr) \, m(\di s) = 0 \nonumber \\
\text{and} \quad
\int _{\uprho (\lambda)}^{\varrho (\lambda)} \bigl( &
K(s) + h(s) - \lambda \bigr) \, m(\di s) > 0 ,
\nonumber
\end{align}
which imply that there exists a unique point $\alpha
(\lambda) \in \mbox{} \bigl] 0 , \uprho (\lambda) \bigr[$
such that
\be
\int _{\alpha (\lambda)}^{\varrho (\lambda)} \bigl(
K(s) + h(s) - \lambda \bigr) \, m(\di s) = 0 .
\ee
For this choice of $\alpha (\lambda)$, we can see
that the function $w_\lambda'$ defined by
(\ref{w-lam-prop1}) and (\ref{w-lam-prop2}) is
indeed $C^1$.
In particular, the limits in (\ref{wlam->w}) all
hold true.
Furthermore,
\ben
\int _{\alpha (\lambda)}^x \bigl( K(s) + h(s) - \lambda
\bigr) \, m(\di s) < 0 \quad \text{for all } x \in \mbox{}
\bigl] \alpha (\lambda), \varrho (\lambda) \bigr[
. \label{wlam-ineq3}
\een

The inequality (\ref{w-lam-ineq1}) follows
immediately from (\ref{wlam-ineq3}).
On the other hand, it is straightforward to
check that the function $w_\lambda'$ defined
by (\ref{w-lam-prop1}) and (\ref{w-lam-prop2})
satisfies the equality in (\ref{w-lam-ineq2}),
while the inequality in (\ref{w-lam-ineq2}) is
equivalent to
\be
K(x) + h(x) - \lambda < 0 , \quad \text{for }
x \in \mbox{} \bigl] 0, \alpha (\lambda) \bigr[
\mbox{} \cup  \mbox{} \bigl] \varrho (\lambda)
, \infty \bigr[ ,
\ee
which is true thanks to the inequalities
(\ref{Kbhl-rho-prop}), (\ref{Kbhl-uprho-prop})
and the fact that $\alpha (\lambda) \in \mbox{}
\bigl] 0 , \uprho (\lambda) \bigr[$.

Finally, (\ref{w'-bound}) follows immediately from
the continuity of $w_\lambda'$ and the boundedness
of $k$.
\end{proof}

\begin{rem} \label{rem:AH19}
The model studied by Alvarez and Hening~\cite{AH22}
is the special case that arises when $h = 0$ and
$k = 1$.
In this case, the identity (\ref{Kk-rel}) in
Remark~\ref{rem:Kk-rel} impies that
\ben
\int _0^x K(s) \, m(\di s) = \int _0^x b(s) \, m(\di s)
= \frac{1}{p'(x)} . \label{Kb-rel}
\een
In view of this identity, we can see that the system
of equations in (\ref{lam-beta-range}), which
determines $(\beta^\star, \lambda^\star)$, and
(\ref{w*-prop}) reduce to
\be
\lambda = b(\beta) , \quad
\lambda = \frac{1}{p' (\beta) m \bigl( ]0 ,\beta[ \bigr)}
\quad \text{and} \quad
w' (x) = \begin{cases} \lambda^\star p' (x) m \bigl( ]0 ,x[ \bigr)
, & \text{for } x \in \mbox{} ]0, \beta^\star[ , \\ 1 , & \text{for }
x \geq \beta^\star , \end{cases}
\ee
which are precisely the expressions (8) and (9) in
Alvarez and Hening~\cite{AH22}.
\end{rem}

\begin{rem} \label{rem:why(III)}
Consider the function $w'$ defined by (\ref{w*-prop})
and suppose that $X$ is as in Example~\ref{ex:VPdif}.
Using L'H\^{o}pital's formula, we calculate
\begin{align}
\lim _{x \downarrow 0} \bigl( x w' (x) \bigr)
& = - \lim _{x \downarrow 0} \frac{\frac{\di}{\di x}
\bigl( x \int _0^x \bigl( K(s) + h(s) - \lambda^\star
\bigr) \, m (\di s) \bigr)}
{\frac{\di}{\di x} \bigl( 1 / p'(x) \bigr)} \nonumber \\
& \geq - \lim _{x \downarrow 0}
\frac{K(x) + h(x) - \lambda^\star}
{\upkappa (\upgamma -x)} =
\frac{\lambda^\star - K(0) - h(0)}{\upkappa \upgamma}
> 0 . \nonumber
\end{align}
It follows that, in the context of
Example~\ref{ex:VPdif}, $\lim _{x \downarrow 0}
w(x) = -\infty$.
\end{rem}

\section{The solution to the ergodic harvesting problem}
\label{sec:erg-sol}

\begin{theorem} \label{VT-th}
Consider the ergodic control problems formulated in
Section~\ref{sec:pr-form}, and let
$(\beta^\star , \lambda^\star)$ be as in
Proposition~\ref{SolPb-lem}.
Given any $x>0$, the following statements hold true:
\smallskip

\noindent {\rm (I)}
$\JE _x (\zeta) \leq \lambda^\star$ and
$\JP _x (\zeta) \leq \lambda^\star$ for all admissible
harvesting strategies $\zeta \in \acal$.
\smallskip

\noindent {\rm (II)}
If $\zeta^\star \in \acal$ is the harvesting strategy that has
a jump of size $\Delta \zeta _0^\star = (x-\beta^\star)^+$
at time 0 and then reflects the state process $X^\star$ at
the level $\beta^\star$ in the negative direction, then
\begin{align}
\JE _x (\zeta^\star) = \lim _{T \uparrow \infty} \frac{1}{T}
\EXP \biggl[ \int _0^T h (X_t^\star) \, \di t + \int _0^T k(X_t^\star)
\circ \di \zeta_t^\star \biggr] & = \lambda^\star \nonumber \\
\text{and} \quad
\JP _x (\zeta^\star) = \lim _{T \uparrow \infty} \frac{1}{T}
\biggl( \int _0^T h (X_t^\star) \, \di t + \int _0^T k(X_t^\star)
\circ \di \zeta_t^\star \biggr)
& = \lambda^\star . \nonumber
\end{align}
\end{theorem}

\begin{proof}
Fix any initial state $x > 0$, let $\zeta \in \acal$ be any
admissible harvesting strategy and let $X$ be the
associated solution to the SDE (\ref{SDE}).
Also, consider the function $w_\lambda$ defined
by (\ref{w-lam-prop1}) and (\ref{w-lam-prop2}) for
$\lambda \in \mbox{} ]\lambda^\star ,
\overline{\lambda}[$.
Using It\^{o}'s formula, we calculate
\begin{align}
w_\lambda (X_T^\zeta) = \mbox{} & w_\lambda (x)
+ \int _0^T \biggl( \frac{1}{2} \sigma^2 (X_t^\zeta)
w_\lambda'' (X_t^\zeta) + b(X_t^\zeta) w_\lambda'
(X_t^\zeta) \biggr) \, \di t - \int _{[0,T]} w_\lambda'
(X_{t-}^\zeta) \, \di \zeta_t \nonumber \\
& + \sum _{0 \leq t \leq T} \bigl( w_\lambda (X_t^\zeta)
- w_\lambda (X_{t-}^\zeta) - w_\lambda' (X_{t-}^\zeta)
\, \Delta X_t^\zeta \bigr) + M_T^\zeta , \nonumber
\end{align}
where
\be
M_T^{\lambda,\zeta} = \int _0^T \sigma (X_t^\zeta)
w_\lambda' (X_t^\zeta) \, \di W_t .
\ee
Since $\Delta X_t^\zeta = X_t^\zeta - X_{t-}^\zeta
= - \Delta \zeta_t \leq 0$ and
\be
w_\lambda (X_t^\zeta) - w_\lambda (X_{t-}^\zeta) +
\int _0^{\Delta \zeta_t} k(X_{t-}^\zeta - u) \, \di u
= \int_{X_t^\zeta}^{X_{t-}^\zeta}
\bigl( k(u) - w_\lambda' (u) \bigr) \, \di u ,
\ee
it follows that
\begin{align}
\int _0^T h(X_t^\zeta) & \, \di t + \int _0^T k(X_t^\zeta)
\circ \di \zeta_t \nonumber \\
= \mbox{} & \lambda T + w_\lambda (x) - w_\lambda
(X_T^\zeta) \nonumber \\
& + \int _0^T \biggl( \frac{1}{2} \sigma ^2
(X_t^\zeta) w_\lambda'' (X_t^\zeta) + b(X_t^\zeta)
w_\lambda' (X_t^\zeta) + h(X_t^\zeta) - \lambda
\biggr) \di t
\nonumber \\
& + \int _0^T \bigl( k(X_t^\zeta) - w_\lambda'
(X_t^\zeta) \bigr) \, \di \zeta_t^{\mathrm c} + \sum
_{0 \leq t \leq T} \int _{X_t^\zeta}^{X_{t-}^\zeta}
\bigl( k(u) - w_\lambda' (u) \bigr) \, \di u + M_T^{\lambda,\zeta}
. \nonumber
\end{align}
Since $\zeta^c$ is an increasing process, $X_t^\zeta
< X_{t-}^\zeta$ and the pair $(w_{\lambda^\star} ,
\lambda^\star)$ (resp., $(w_\lambda , \lambda)$)
satisfies the HJB equation (\ref{HJBe}) (resp., the
inequalities (\ref{w-lam-ineq1}) and (\ref{w-lam-ineq2})),
we can see that
\ben
\int _0^T h(X_t^\zeta) \, \di t + \int _0^T k(X_t^\zeta)
\circ \di \zeta_t \leq \lambda T + w_\lambda
(x) - w_\lambda (X_T^\zeta) + M_T^{\lambda,\zeta}
. \label{wXT}
\een

\underline{\em Proof of the inequality
$\JE _x (\zeta) \leq \lambda^\star$\/.}
Fix any $\lambda \in \mbox{} ]\lambda^\star ,
\overline{\lambda}[$ and let $(\tau_n)$ be a sequence
of localising times for the corresponding local
martingale $M^{\lambda,\zeta}$.
Recalling the assumptions that $h$ is bounded from below
and $k$ is positive, as well as the facts that $\zeta$ is an
increasing process and $w_\lambda$ is bounded from below
(see (\ref{w'-bound}) in Proposition~\ref{SolPb-lem}.(III)),
we take expectations in (\ref{wXT}) and we use the
monotone and the dominated convergence theorems to
calculate
\begin{align}
\frac{1}{T} \EXP \biggl[ \int _0^T h(X_t^\zeta) & \, \di t
+ \int _0^T k(X_t^\zeta) \circ \di \zeta_t \biggr]
\nonumber \\
& = \frac{1}{T} \lim _{n \uparrow \infty} \EXP \biggl[
\int _0^{\tau_n \wedge T} h(X_t^\zeta) \, \di t +
\int _0^{\tau_n \wedge T} k(X_t^\zeta)
\circ \di \zeta_t \biggr] \nonumber \\
& \leq \frac{1}{T} \lim _{n \uparrow \infty} \EXP \Bigl[
\lambda (\tau_n \wedge T) + w_\lambda (x) + w_\lambda^-
(X_{\tau_n \wedge T}^\zeta) \Bigr] \nonumber \\
& = \lambda + \frac{w_\lambda (x)}{T} + \frac{1}{T}
\EXP \Bigl[ w_\lambda^- (X_T^\zeta) \Bigr] , \nonumber
\end{align}
where $w_\lambda^- (x) = - \min \bigl\{ w_\lambda(x)
, 0 \bigr\}$.
Using the fact that $w_\lambda^-$ is bounded once
again, we can pass to the limit as $T \uparrow \infty$
to obtain the inequality $\JE _x (\zeta) \leq \lambda$,
which implies the required inequality
$\JE _x (\zeta) \leq \lambda^\star$ by passing to the
limit as $\lambda \downarrow \lambda^\star$.
\smallskip

\underline{\em Proof of the inequality
$\JP _x (\zeta) \leq \lambda^\star$\/.}
Making a slight modification of the proof of the comparison
Theorem~V.43 in Rogers and Williams~\cite{RW},
we can show that $X_t^\zeta \leq X_t$ for all
$t \geq 0$, $\PR$-a.s., where $X$ is
the solution to the SDE (\ref{SDE-0}).
In view of this observation, we can see that, given
any $\lambda \in \mbox{} ]\lambda^\star ,
\overline{\lambda}[$,
\begin{align}
\bigl\langle M^{\lambda,\zeta} \bigr\rangle _T
& = \int _0^T \Bigl( \sigma (X_t^\zeta) w_\lambda'
(X_t^\zeta) \Bigr) ^2 \, \di t \nonumber \\
& \leq C_1 C_2^2 \int _0^T \Bigl( 1 + \bigl(
X_t^\zeta \bigr)^\eta \Bigr) \, \di t \leq C_1 C_2^2
\biggl( T + \int _0^T X_t^\eta \, \di t \biggr) ,
\nonumber
\end{align}
where $C_1$, $\eta$ and $C_2 = C_2 (\lambda)$ are
the constants in (\ref{sigma-bounds}) and
(\ref{w'-bound}).
Furthermore, the ergodic Theorem~V.53 in
Rogers and Williams~\cite{RW} implies that
\begin{align}
\begin{split}
\limsup _{T \uparrow \infty} \frac{\bigl\langle
M^{\lambda,\zeta} \bigr\rangle _T}{T} & \leq
C_1 C_2^2 \biggl( 1 + \lim _{T \uparrow \infty}
\frac{1}{T} \int _0^T X_t^\eta \, \di t \biggr) \\
& = C_1 C_2^2 \biggl( 1 + \frac{1}{m \bigl( ]0,
\infty[ \bigr)} \int _0^\infty s^\eta \, m (\di s)
\biggr) =: C_3 < \infty ,
\end{split}
\label{<M>}
\end{align}
with the second inequality following thanks to
Assumption~\ref{A12+}.

The Dambis, Dubins and Schwarz theorem
(e.g., see Revuz and Yor~\cite[Theorem V.1.7]{RY}) asserts
that there exists a standard Brownian motion $B$, which
may be defined on a possible enlargement of the probability
space $(\Omega, \fcal, \PR)$, such that $M^{\lambda,\zeta}
= B_{\langle M^{\lambda,\zeta} \rangle}$.
Using this representation, (\ref{<M>}) and the fact that
$\lim _{T \uparrow \infty} B_T / T = 0$, we can see
that
\be
\lim _{T \uparrow \infty}
\frac{\bigl| M_T^{\lambda,\zeta} \bigr|}{T} {\bf 1}
_{\{ \langle M^{\lambda,\zeta} \rangle _\infty = \infty \}}
\leq C_3 \lim _{T \uparrow \infty}
\frac{|B_{\langle M^{\lambda,\zeta} \rangle_T}|}
{\bigl\langle M^{\lambda,\zeta} \bigr\rangle _T} {\bf 1}
_{\{ \langle M^{\lambda,\zeta} \rangle _\infty = \infty \}}
= 0 .
\ee
On the other hand,
\be
\lim _{T \uparrow \infty}
\frac{\bigl| M_T^{\lambda,\zeta} \bigr|}{T} {\bf 1}
_{\{ \langle M^{\lambda,\zeta} \rangle _\infty < \infty \}}
= 0
\ee
because $M^{\lambda,\zeta}$ converges in $\bbr$
on the event $\{ \bigl\langle M^{\lambda,\zeta}
\bigr\rangle _\infty < \infty \}$.
In view of these results, we can pass to the limit as
$T \uparrow \infty$ in (\ref{wXT}) to obtain
\be
\JP _x (\zeta) \leq \lim _{T \uparrow \infty} \biggl(
\lambda + \frac{w_\lambda (x)}{T} +
\frac{w_\lambda^- (X_T^\zeta)}{T}
+ \frac{M_T^{\lambda,\zeta}}{T} \biggr) = \lambda .
\ee
The inequality $\JP _x (\zeta) \leq \lambda^\star$
now follows by passing to the limit as $\lambda
\downarrow \lambda^\star$.
\smallskip

\underline{\em Proof of (II)\/.}
Let the harvesting strategy $\zeta^\star \in \acal$ be as
in the statement of the theorem: such a strategy indeed
exists (see Tanaka~\cite[Theorem~4.1]{Tan}).
If we define
\be
N_T = \int _0^T \sigma (X_t^\star) \, \di W_t ,
\ee
then
$\langle N \rangle _T / T \leq \max _{s \in [0,\beta^\star]}
\sigma (s) < \infty$.
Therefore, $N$ is a square integrable martingale
and $\EXP \bigl[ N_T \bigr] = 0$ for all $T>0$.
Furthermore, reasoning as above, we can see that
$\lim _{T \uparrow \infty} N_T / T = 0$.
In view of these observations, the expression
\ben
\frac{\zeta_T^\star}{T} = \frac{x}{T} - \frac{X_t^\star}{T}
+ \frac{1}{T} \int _0^T b(X_t^\star) \, \di t  + \frac{N_T}{T} ,
\een
and the fact that, beyond its possible initial jump,
$\zeta^\star$ increases on the set $\{ X_t^\star =
\beta^\star \}$, we can see that
\begin{align}
\JE _x (\zeta^\star) & = \lim _{T \uparrow \infty} \frac{1}{T}
\EXP \biggl[ \int _0^T h (X_t^\star) \, \di t + k(\beta^\star)
\zeta_T^\star \biggr] \nonumber \\
& = \lim _{T \uparrow \infty} \frac{1}{T}
\EXP \biggl[ \int _0^T \Bigl( h (X_t^\star) + k(\beta^\star)
b (X_t^\star) \Bigr) \, \di t \biggr] \nonumber
\end{align}
and
\be
\JP _x (\zeta^\star) = \lim _{T \uparrow \infty} \frac{1}{T}
\biggl( \int _0^T h (X_t^\star) \, \di t + k(\beta^\star)
\zeta_T^\star \biggr) = \lim _{T \uparrow \infty} \frac{1}{T}
\int _0^T \Bigl( h (X_t^\star) + k(\beta^\star) b (X_t^\star)
\Bigr) \, \di t .
\ee
These expressions and standard ergodic theorems
(e.g., see Borodin and Salminen~\cite[Section~II.6]{BS}
and Rogers and Williams~\cite[Theorem~V.53]{RW})
imply that
\be
\JE _x (\zeta^\star) = \JP _x (\zeta^\star) =
\frac{1}{m \bigl( ]0, \beta^\star[ \bigr)}
\int _0^{\beta^\star} \bigl( h(s) + k(\beta^\star) b(s)
\bigr) \, m(\di s) .
\ee
Combining these observations with the identities
\be
k(\beta^\star) \int _0^{\beta^\star} b(s) \, m(\di s)
\stackrel{(\ref{Kb-rel})}{=}
\frac{k(\beta^\star)}{p' (\beta^\star)}
\stackrel{(\ref{Kk-rel})}{=} \int _0^{\beta^\star}
K(s) \, m(\di s) ,
\ee
we obtain
 \be
\JE _x (\zeta^\star) = \JP _x (\zeta^\star) =
\frac{1}{m \bigl( ]0, \beta^\star[ \bigr)}
\int _0^{\beta^\star} \bigl( K(s) + h(s) \bigr) \, m(\di s)
\stackrel{(\ref{lambda*location})}{=} \lambda^\star .
\ee
\end{proof}

\section{Abelian limits}
\label{sec:Abelian}

In this section, we allow for the discounting rate function
$r$ to depend on a parameter $\iota > 0$ and we establish
the convergence of the solution to the discounted
control problem to the one of the ergodic control
problems in an Abelian sense.
In particular, we make the following assumption, which
is the same as Assumption~\ref{A4d}.(i) for each
individual $\iota > 0$.

\begin{ass} \label{ass:Ab}
The discounting rate function $(x,\iota) \mapsto
r(x; \iota)$ is continuous.
Also, given any $\iota > 0$, the function
$r(\cdot; \iota)$ is $C^1$ and such that
\ben
\underline{r} (\iota) \leq r(x; \iota) \leq \overline{r}
(\iota) \quad \text{for all } x \geq 0 , \label{rxp-bounds}
\een
for some $\underline{r} (\iota)$ and $\overline{r}
(\iota)$ such that 
\ben
0 < \underline{r} (\iota) < \overline{r} (\iota)
< \infty \text{ for all } \iota > 0 , \quad \lim
_{\iota \downarrow 0} \frac{\overline{r} (\iota)}
{\underline{r} (\iota)} = 1 \quad \text{and} \quad
\lim _{\iota \downarrow 0} \overline{r} (\iota)
= 0 . \label{rp-props}
\een
\end{ass}
The dependence of $r$ on the parameter $\iota$ implies
that the functions $\Theta$, $\varphi$, $\psi$ and
$R_h$ that we have considered in our analysis also
depend on $\iota$.
Throughout this section, we will make such
dependences explicit for clarity of the arguments.

The functions $\varphi$ and $\psi$ introduced at the
beginning of Section~\ref{sec:disc-sol} are unique
up to a multiplicative constant.
In this section, we assume that they have been
scaled so that
\ben
\varphi (1; \iota) = 1 \quad \text{and} \quad
\psi (1; \iota) = 1 \quad \text{for all } \iota > 0 ,
\label{phi-psi-scaled}
\een
without loss of generality.

\begin{lemma}
In the presence of Assumptions~\ref{A1} and~\ref{A2},
the scaled as in (\ref{phi-psi-scaled}) functions
$(x, \iota) \mapsto \varphi (x; \iota)$ and
$(x, \iota) \mapsto \psi (x; \iota)$ are continuous,
\begin{gather}
\lim _{\iota \downarrow 0} \frac{\varphi (x; \iota)}
{\varphi (y; \iota)} = \lim _{\iota \downarrow 0}
\frac{\psi (x; \iota)}{\psi (y; \iota)} = 1 \quad \text{for all }
x,y > 0 \label{phi-psi-plims} \\
\text{and} \quad
\lim _{\iota \downarrow 0} \frac{\psi' (x; \iota)}
{r(y; \iota) \psi (x; \iota)} = p'(x) m \bigl( ]0,x[ \bigr)
\quad \text{for all } x,y > 0 . \label{psi'/psi-plim}
\end{gather}
In the presence of the assumptions we have made
in Section~\ref{sec:pr-form} and
Assumption~\ref{ass:Ab},
\ben
\lim _{\iota \downarrow 0} r(y; \iota)
\biggl( R_h (x; \iota) - \frac{R_h' (\beta; \iota)}
{\psi' (\beta; \iota)} \psi (x; \iota) \biggr)
= \frac{1} {m \bigl( ]0, \beta[ \bigr)}
\int _0^\beta h(s) \, m (\di s) \label{reflRh-Ablim}
\een
for all $\beta > 0$, $x \in \mbox{} ]0, \beta]$ and
$y > 0$.
\end{lemma}

\begin{proof}
The continuity of the functions $\varphi$ and $\psi$,
as well as (\ref{phi-psi-plims}), follow immediately from
(\ref{phi-psi-defn}) and the dominated convergence
theorem.
In turn, (\ref{psi'/psi-plim}) follows from the definition
(\ref{Psi}), the identity (\ref{psi/p'-int}), the limit
(\ref{phi-psi-plims}) and Assumption~\ref{ass:Ab},
which imply that
\be
\lim _{\iota \downarrow 0} \frac{\psi' (x; \iota)}
{r(y; \iota) \psi (x; \iota)} = p'(x) \lim
_{\iota \downarrow 0} \int _0^x
\frac{r(s; \iota)}{r(y; \iota)} \frac{\varphi (s; \iota)}
{\varphi (x; \iota)} \, m (\di s) = p'(x) m
\bigl( ]0,x[ \bigr) .
\ee
Using the definitions (\ref{Rh})--(\ref{Psi}), we
can see that
\begin{align}
R_h & (x; \iota) - \frac{R_h' (\beta; \iota)}
{\psi' (\beta; \iota)} \psi (x; \iota) \nonumber \\
& = \frac{1}{C} \widetilde{\varphi} _\beta (x; \iota)
\int _0^x h(s) \psi (s; \iota) \, m (\di s) + \frac{1}{C}
\psi (x; \iota) \int _x^\beta h(s) \widetilde{\varphi} _\beta
(s; \iota) \, m(\di s) , \nonumber
\end{align}
where
\be
\widetilde{\varphi} _\beta (x; \iota) = \varphi (x; \iota)
- \frac{\varphi' (\beta; \iota)}{\psi' (\beta; \iota)} \psi
(x; \iota) .
\ee
In view of the observation that
\begin{align}
\frac{\widetilde{\varphi} _\beta (x; \iota)}
{\widetilde{\varphi} _\beta (\beta; \iota)}
= \mbox{} & 1 + \biggl( \frac{\varphi (x; \iota)}
{\varphi (\beta; \iota)} - 1 \biggr)
\frac{\varphi (\beta; \iota) \psi' (\beta; \iota)}
{\varphi (\beta; \iota) \psi' (\beta; \iota) -
\varphi' (\beta; \iota) \psi (\beta; \iota)}
\nonumber \\
& + \biggl( \frac{\psi (x; \iota)}{\psi (\beta; \iota)}
- 1 \biggr) \frac{- \varphi' (\beta; \iota) \psi
(\beta; \iota)}{\varphi (\beta; \iota) \psi'
(\beta; \iota) - \varphi' (\beta; \iota) \psi
(\beta; \iota)} \nonumber
\end{align}
and the fact that the two long fractions on the right-hand
side of this expression take values in $]0,1[$,
we can see that $\lim _{\iota \downarrow 0}
\widetilde{\varphi} _\beta (x; \iota) / \widetilde{\varphi}
_\beta (\beta; \iota) = 1$, thanks to
(\ref{phi-psi-plims}).
On the other hand, we can use the probabilistic expression
in (\ref{Rh}) to obtain
\be
R_{r (\cdot; \iota)} (x) = \EXP_x \biggl[ \int _0^\infty
\exp \biggl( - \int _0^t r(X_u ; \iota) \, \di u \biggr)
r(X_t ; \iota) \, \di t \biggr] = 1 \quad \text{for all }
x, \iota > 0 .
\ee
Combining these observations with (\ref{phi-psi-plims})
and the fact that
\ben
\lim _{\iota \downarrow 0} \frac{r(x; \iota)}{r(y; \iota)}
= 1 \quad \text{for all } x,y > 0 , \label{r-plims}
\een
which follows from Assumption~\ref{ass:Ab}, and
using the dominated convergence theorem, we obtain
\begin{align}
\lim _{\iota \downarrow 0} r(y; \iota) \biggl( & R_h
(x; \iota) - \frac{R_h' (\beta; \iota)}{\psi' (\beta; \iota)}
\psi (x; \iota) \biggr) \nonumber \\
& = \lim _{\iota \downarrow 0} r(y; \iota)
\frac{R_h (x; \iota) - \frac{R_h' (\beta; \iota)}
{\psi' (\beta; \iota)} \psi (x; \iota)}
{R_{r (\cdot; \iota)} (x;p) -
\frac{R_{r (\cdot; \iota)}' (\beta; \iota)}
{\psi' (\beta; \iota)} \psi (x; \iota)} \nonumber \\
& = \lim _{\iota \downarrow 0}
\frac{\int _0^x h(s) \frac{\psi (s; \iota)}{\psi (x; \iota)}
\, m (\di s) + \int _x^\beta h(s)
\frac{\widetilde{\varphi} _\beta (s; \iota)}
{\widetilde{\varphi} _\beta (x; \iota)} \, m(\di s)}
{\int _0^x \frac{r(s; \iota)}{r(y; \iota)}
\frac{\psi (s; \iota)}{\psi (x; \iota)} \, m (\di s)
+ \int _x^\beta \frac{r(s; \iota)}{r(y; \iota)}
\frac{\widetilde{\varphi} _\beta (s; \iota)}
{\widetilde{\varphi} _\beta (x; \iota)} \, m(\di s)}
= \frac{\int _0^\beta h(s) \, m (\di s)}
{m \bigl( ]0, \beta[ \bigr)} , \nonumber
\end{align}
namely, (\ref{reflRh-Ablim}).
\end{proof}

\begin{theorem}
Consider the control problems formulated in
Section~\ref{sec:pr-form} and suppose that
Assumption~\ref{ass:Ab} also holds.
If $\beta^\star (\iota)$, $w (\cdot ; \iota)$ are
as in Lemma~\ref{lem-existGamma} and
$(\beta^\star , \lambda^\star)$, $w$ are as in
Proposition~\ref{SolPb-lem}, then
\ben
\lim _{\iota \downarrow 0} \beta^\star (\iota)
= \beta^\star , \quad \lim _{\iota \downarrow 0}
r(y; \iota) w(x; \iota) = \lambda^\star
\quad \text{and} \quad
\lim _{\iota \downarrow 0} w'(x; \iota)
= w'(x) \label{Ab-lims}
\een
for all $x,y > 0$.
\end{theorem}

\begin{proof}
In view of the definition (\ref{Psi}), the equation
(\ref{FA}) that $\beta^\star (\iota) > 0$ satisfies takes the
form
\begin{align}
\int _0^{\beta^\star (\iota)} \bigl( \Theta (s; \iota)
& + h(s) \bigr) \frac{\psi (s; \iota)}
{\psi \bigl( \beta^\star (\iota) ; \iota \bigr)}
\, m (\di s) \nonumber \\
& = \bigl( \Theta \bigl( \beta^\star (\iota) ; \iota \bigr)
+ h \bigl( \beta^\star (\iota) \bigr) \bigr) \int
_0^{\beta^\star (\iota)}
\frac{r(s; \iota)}{r \bigl( \beta^\star (\iota) ; \iota \bigr)}
\frac{\psi (s; \iota)}
{\psi \bigl( \beta^\star (\iota) ; \iota \bigr)}
\, m (\di s) . \nonumber
\end{align}
The functions $r$, $\Theta$, $h$ and $\psi$ are all
continuous, while $\lim _{\iota \downarrow 0} \Theta
(x; \iota) = K(x)$ (see the definitions (\ref{Theta}) and
(\ref{K}) of $\Theta$ and $K$).
Therefore, we can use (\ref{phi-psi-plims}), (\ref{r-plims})
and the dominated convergence theorem  to come to
the conclusion that the limit
$\beta^\star (0) = \lim _{\iota \downarrow 0} \beta^\star
(\iota)$ exists and satisfies the equation
\be
\int _0^{\beta^\star (0)} \bigl( K(s) + h(s) \bigr)
\, m (\di s)
= \bigl( K \bigl( \beta^\star (0) \bigr) +
h \bigl( \beta^\star (0) \bigr) \bigr)
m \bigl( \bigl] 0, \beta^\star (0) \bigr[ \bigr) .
\ee
It follows that the first limit in (\ref{Ab-lims}) holds true
because this is the equation that $\beta^\star$
uniquely satisfies (see (\ref{lambda*location})).

The second limit in (\ref{Ab-lims}) follows
immediately from the first expression for
$w (\cdot; \iota)$ in (\ref{OS3}), Assumption~\ref{ass:Ab}
and (\ref{reflRh-Ablim}) with $\Theta (\cdot; \iota) + h$
in place of $h$.
Finally, we use (\ref{C}) and (\ref{phi-psi-plims})
to note that
\begin{align}
\lim _{\iota \downarrow 0} \biggl( R _{\Theta +h}'
- \frac{\psi'}{\psi} R _{\Theta +h} \biggr) (x; \iota)
& = - p' (x) \lim _{\iota \downarrow 0} \int _0^x \bigl(
\Theta (s; \iota) + h(s) \bigr) \frac{\psi (s; \iota)}
{\psi (x ; \iota)} \, m (\di s) \nonumber \\
& = - p' (x) \int _0^x \bigl( K(s) + h(s) \bigr) \, m (\di s) .
\nonumber
\end{align}
In light of this limit, the fact that $\lim
_{\iota \downarrow 0} \beta^\star (\iota) = \beta^\star$,
(\ref{psi'/psi-plim}) and (\ref{reflRh-Ablim}) with
$\Theta (\cdot; \iota) + h$ in place of $h$, we
can see that
\begin{align}
\lim _{\iota \downarrow 0} \biggl( & R _{\Theta +h}'
(x; \iota) -
\frac{R_{\Theta +h}' \bigl( \beta^\star (\iota); \iota \bigr)}
{\psi' \bigl( \beta^\star (\iota); \iota \bigr)} \psi' (x; \iota)
\biggr) \nonumber \\
= \mbox{} & \lim _{\iota \downarrow 0} \biggl(
R _{\Theta +h}' - \frac{\psi'}{\psi} R _{\Theta +h}
\biggr) (x; \iota) \nonumber \\
& + \lim _{\iota \downarrow 0} r(y; \iota)
\biggl( R _{\Theta +h} (x; \iota) -
\frac{R _{\Theta +h}' \bigl( \beta^\star (\iota); \iota \bigr)}
{\psi' \bigl( \beta^\star (\iota); \iota \bigr)} \psi
(x; \iota) \biggr) \lim _{\iota \downarrow 0}
\frac{\psi' (x; \iota)}{r(y; \iota) \psi (x; \iota)}
\nonumber \\
= \mbox{} & p'(x) \Biggl(
\frac{m \bigl( ]0, x[ \bigr)}{m \bigl( ]0, \beta^\star[ \bigr)}
\int _0^{\beta^\star} \bigl( K(s) + h(s) \bigr) \, m
(\di s) - \int _0^x \bigl( K(s) + h(s) \bigr) \, m (\di s)
\Biggr) . \nonumber
\end{align}
The third limit in (\ref{Ab-lims}) follows from this
result, the first expression for $w(\cdot; \iota)$
in (\ref{OS3}) and the second expression for
$w'$ in (\ref{w'1}).
\end{proof}

\section*{Acknowledgment}

We are grateful to three anonymous referees whose
constructive comments led to a substantial
enhancement of the paper.



\begin{thebibliography}{20}
%
\bibitem{ALV00}
{\sc L. H. R. Alvarez},
{\em On the option interpretation of rational harvesting planning},
J. Math. Biol.,
40 (2000), pp. 383--405.
%
\bibitem{ALV01}
{\sc L. H. R. Alvarez},
{\em Singular stochastic control, linear diffusions, and
optimal stopping: A class of solvable problems},
SIAM J. Control Optim.,
39 (2001), pp. 1697--1710.
%
\bibitem{AH22}
{\sc L. H. R. Alvarez and A.\,Hening},
{\em Optimal sustainable harvesting of populations in
random environments},
Stochastic Process. Appl.,
150 (2022), pp. 678--698.
%
\bibitem{ALO}
{\sc L. H. R. Alvarez, E. Lungu and B. {\O}ksendal},
{\em Optimal multi-dimensional stochastic harvesting
with density-dependent prices},
Afr. Mat.,
27 (2016), pp. 427--442.
%
\bibitem{AS}
{\sc L. H. R. Alvarez and L. A. Shepp},
{\em Optimal harvesting of stochastically fluctuating
populations},
J. Math. Biol.,
37 (1998), pp. 155--177.
%
\bibitem{BS}
{\sc A. N. Borodin and P. Salminen},
{\em Handbook of Brownian Motion - Facts and Formulae\/},
Birkh\"{a}user, 2002.
%
\bibitem{CDF23}
{\sc H. Cao, J. Dianetti and G. Ferrari},
{\em Stationary discounted and ergodic mean field
games of singular control},
Math. Oper. Res.,
48 (2023), pp. 1871--1898.
%
\bibitem{CHS22}
{\sc A. Cohen, A. Hening and C. Sun},
{\em Optimal ergodic harvesting under ambiguity},
SIAM J. Control Optim.,
60 (2022), pp. 1039--1063.
%
\bibitem{FRA}
{\sc N. C. Framstad},
{\em Optimal harvesting of a jump diffusion population
and the effect of jump uncertainty},
SIAM J. Control Optim.,
42 (2003), pp. 1451--1465.
%
\bibitem{GLVS22}
{\sc M. Ga\"{i}gi, V. Ly Vath and S. Scotti},
{\em Optimal harvesting under marine reserves and
uncertain environment},
European J. Oper. Res.,
301 (2022), pp. 1181--1194.
%
\bibitem{HNUW}
{\sc A. Hening, D. H. Nguyen, S. C. Ungureanu
and T. K. Wong},
{\em Asymptotic harvesting of populations in random
environments},
J. Math. Biol.,
78 (2019), pp. 293--329.
%
\bibitem{HTPY}
{\sc A. Hening, K. Q. Tran, T. T. Phan and G. Yin},
{\em Harvesting of interacting stochastic populations},
J. Math. Biol.,
79 (2019), pp. 533--570.
%
\bibitem{H12}
{\sc R. Hynd},
{\em The eigenvalue problem of singular ergodic control},
Comm. Pure Appl. Math.,
65 (2012), pp. 649--682.
%
\bibitem{JZ}
{\sc A. Jack and M. Zervos},
{\em A singular control problem with an expected
and a pathwise ergodic performance criterion},
J. Appl. Math. Stochastic Anal.,
 (2006), article ID 82538, 19 pp.
%
\bibitem{K83}
{\sc I. Karatzas},
{\em A class of singular stochastic control problems},
Adv. in Appl. Probab.,
15 (1983), pp. 225--254.
%
\bibitem{KSh}
{\sc I. Karatzas and S. E. Shreve},
{\em Brownian Motion and Stochastic Calculus},
Springer-Verlag, 1988.
%
\bibitem{KXYZ22}
{\sc K. Kunwai, F. Xi, G. Yin and C. Zhu},
{\em On an ergodic two-sided singular control problem},
Appl. Math. Optim.,
86 (2022), paper no. 26, 34 pp.
%
\bibitem{LZ}
{\sc D. Lamberton and M. Zervos},
{\em On the optimal stopping of a one-dimensional diffusion},
Electron. J. Probab.,
18 (2013), paper no. 34, 49 pp.  
%
\bibitem{LZ24}
{\sc Z. Liu and M. Zervos},
{\em The solution to an impulse control problem
motivated by optimal harvesting},
J. Math. Anal. Appl.,
542 (2025), paper no. 128809, 33 pp.
%
\bibitem{LO97}
{\sc E. Lungu and B. {\O}ksendal},
{\em Optimal harvesting from a population in a
stochastic crowded environment},
Math. Biosci.,
145 (1997), pp. 47--75.
%
\bibitem{LO01}
{\sc E. Lungu and B. {\O}ksendal},
{\em Optimal harvesting from interacting populations
in a stochastic environment},
Bernoulli,
7 (2001), pp. 527--539.
%
\bibitem{M13}
{\sc H. Morimoto},
{\em Optimal harvesting in the stochastic logistic
growth model with finite horizon},
SIAM J. Control Optim.,
51 (2013), pp. 2335--2355.
%
\bibitem{MRT92}
{\sc J. L. Menaldi, M. Robin and M. I. Taksar},
{\em Singular ergodic control for multidimensional
Gaussian processes},
Math. Control Signals Systems,
5 (1992), pp. 93--114.
%
\bibitem{RY}
{\sc D. Revuz and M. Yor},
{\em Continuous Martingales and Brownian Motion},
3rd ed., Springer, 1999.
%
\bibitem{RW}
{\sc L. C. G. Rogers and D. Williams},
{\em Diffusions, Markov Processes and Martingales},
vol. 2, Cambridge University Press, 2000.
%
\bibitem{SSZ}
{\sc Q. Song, R. H. Stockbridge and C. Zhu},
{\em On optimal harvesting problems in random
environments},
SIAM J. Control Optim.,
49 (2011), pp. 859--889.
%
\bibitem{Tan}
{\sc H. Tanaka},
{\em Stochastic differential equations with reflecting
boundary condition in convex regions},
Hiroshima Math. J.,
9 (1979), pp. 163--177.
%
\bibitem{W02}
{\sc A. Weerasinghe},
{\em Stationary stochastic control for It\^o processes},
Adv. in Appl. Probab.,
34 (2002), pp. 128--140.
%
\bibitem{W07}
{\sc A. Weerasinghe},
{\em An Abelian limit approach to a singular
ergodic control problem},
SIAM J. Control Optim.,
46 (2007), pp. 714--737.
%
\end{thebibliography}
\end{document}